\documentclass{elsart}
\usepackage{epsfig}
\usepackage{amssymb,colortbl,euscript}
\usepackage[cp1251]{inputenc}
\usepackage[english]{babel}
\usepackage[T2A]{fontenc}
\usepackage{amsmath}

\begin{document}

\begin{frontmatter}
\title{A Comparison Study  of Two High Accuracy Numerical Methods for a Parabolic System in Air Pollution Modelling}
\author[label1]{I. Dimov},
\author[label2]{J. Kandilarov},
 \author[label1]{V. Todorov},
\author[label2]{L. Vulkov}
\address[label1]{Institute of Information and Communication Technologies,
  BAS, Bulgaria}
\ead{ivdimov@bas.bg}
 \address[label2]{Center of Applied Mathematics and Informatics,
  University of Ruse, Ruse 7017, Bulgaria}
\ead{ukandilarov@uni-ruse.bg}
\ead{venelintodorov@fmi.uni-sofia.bg}
\ead{lvalkov@ru.acad.bg}

\begin{abstract}
We present two approaches for enhancing the accuracy of second order finite difference approximations of two-dimensional semilinear parabolic systems. These are the fourth order compact difference scheme and the fourth order scheme based on Richardson extrapolation. Our interest is concentrated on a system of ten parabolic partial differential  equations in air pollution modeling. We analyze numerical experiments to compare the two approaches with respect to accuracy, computational complexity, non-negativity preserving and etc. Sixth-order approximation based on the fourth-order compact difference scheme combined with Richardson extrapolation is also discussed numerically.
\end{abstract}

\begin{keyword}
air pollution model, semilinear parabolic systems, compact finite difference schemes, Richardson extrapolation.
\end{keyword}

\end{frontmatter}


\section{Introduction}
In many fields of sciences and engineering parabolic equations
 are always used to describe many phenomena, so that the finite-difference method that solves the parabolic equation
is always a focus of concern, see e.g. \cite{CheKin,MarSha,RichtMor,Sam}. In the context of the finite difference discretization, the
standard second-order discretization schemes may need fine griddings to yield approximate solutions of
acceptable accuracy. The resulting large size systems have to be solved, which may consume a lot of memory  space and
CPU cycles even on present generation supercomputers.

One approach to reduce computational cost in very large-scale modelings and simulations is to used higher-order discretization
methods. Other important factor affecting the computational efficacy of a discretized method is to solve the resulting linear and nonlinear systems of algebraic equations.
The higher-order methods usually generate algebraic systems of much smaller size, compared to the lower-order methods.

Because of this and other advantages of high-order methods, there has been growing interest of developing and using highly accurate numerical schemes for solving partial differential equations,
leading to renewed interest in high-order compact difference schemes
\cite{GupManSte,KeyRoo,Rich,SunZha,TurGor,Wang}.

Compact schemes, proposed by Kreiss and Oliger \cite{GusKreOli} use similar stencil, but requires
a scalar tridiagonal or pentadiagonal matrix inversion. In this paper we use another idea
 to obtain high-order
compact schemes, namely, to operate on the differential equations as auxiliary relations in
order to express hight order derivatives in the truncation error \cite{SpoGartime,WangGuoWu}.
More details and discussions on construction of compact difference schemes for convection-diffusion
problems can be found in \cite{KeyRoo,SunZha,WangGuoWu}.

In \cite{KarKur} the air pollution problem, which is the base of the mathematical model of the present paper is stated. A preconditioned iterative solution method for nonlinear parabolic transport system is done. The ingredients of the method are implicit Euler discretization in time and FEM discretization in space, then an outer-inner iteration and preconditioning via an $\ell$-tuple of independent elliptic operators.

Another known approach for increasing the order of accuracy of the finite difference schemes is the use of Richardson extrapolation \cite{MarSha}. Fourth order compact difference scheme for a system of two semilinear toy 1D parabolic equations  is derived in \cite{DKTVHis}.

This article is arranged as follows. In Section 2 we present the two dimensional model problem. In Section 3 the second order central difference scheme (CDS) is presented and the application of the Richardson extrapolation for higher-order approximations is analyzed.
 In Section 4 the fourth-order compact finite difference schemes (CFDS) for general weakly coupled parabolic system of two equations is introduce.
 In Section 5 numerical results and comparisons are presented and analysed. Concluding remarks are included in Section 6.
\section{The Two Dimensional Model Problem of Air Pollution}
The simulation of various processes in chemistry, physics and
engineering  uses models of systems of coupled parabolic problems.
In this work we construct compact high-order finite difference
schemes for semilinear parabolic systems and propose fast algorithms for
solution of the nonlinear algebraic equations. Problems of air
pollution transport with coupling in the nonlinear reactions terms
are of our main consideration, namely,
\begin{equation}\label{2Dsystem1}
 \frac{\partial u_l}{\partial t}
 - K \triangle u_{l} +  \mathbf{b}_{l}  \nabla u_{l}
= R_{l} (x,y, u_{1} , \dots , u_{L} ), \quad (x,y,t)\in  \Omega \times (0,T],
\end{equation}
\begin{equation}\label{2Dsystem2}
\mathbf{u}=0, \quad (x,y,t) \in \partial \Omega \times  (0,T],
\end{equation}
\begin{equation}\label{2Dsystem3}
\mathbf{u}=\mathbf{u}_0(x,y), \quad  \quad    (x,y) \in \Omega,
\end{equation}
\noindent where $\mathbf{u}=(u_1,u_2,...,u_L)$, $u_l=u_l(x,y,t)$, $l=1,...,L$ are the concentrations
of $L$ chemical species (pollutants) and $K>0$ is the diffusion
coefficient and $\Omega \in R^2$ is a bounded domain. The assumption regarding  constant $K: =K_{x} =K_{y} $ is not a
restriction for developing our numerical approach. This just
corresponds to the physical model described in
\cite{ DimZla,GeoZla,KarKur}.

The main goal of the paper is the application and numerical illustration of above-mentioned
difference approximations to the following real-life parabolic transport system described in \cite{GeoZla}.
Following \cite{GeoZla,KarKur,ZlaDim} the advection part in (\ref{2Dsystem1}) may be presented in the following form:
\begin{equation*}
    \mathbf{b_{l}} . \nabla u_{l}=\mu(y-y_c)\frac{\partial u_{l}}{\partial
    x}+\mu(x_c-x)\frac{\partial u_{l}}{\partial y},
\end{equation*}
where $x\in (0,X)$, $y\in (0,Y)$, $x_c=X/2$, $y_c=Y/2$. The
nonlinear chemical part of the model is (see \cite{KarKur}):
\begin{eqnarray}\label{starn}
  R_1(u_1,...,u_{10}) &=& k_5u_2-(k_6u_5+k_4u_7+k_3u_8)u_1,  \nonumber\\
  R_2(u_1,...,u_{10}) &=& (k_6u_5+k_4u_7+k_3u_8)u_1-(k_5+k_9u_9)u_2, \nonumber \\
  R_3(u_1,...,u_{10}) &=& -k_1u_3u_9,  \nonumber\\
  R_4(u_1,...,u_{10}) &=& 2 k_1u_3u_9+k_3u_1u_8-k_2u_4, \nonumber\\
  R_5(u_1,...,u_{10}) &=& k_2u_5 \\
  R_6(u_1,...,u_{10}) &=& k_9u_2u_9,  \nonumber\\
  R_7(u_1,...,u_{10}) &=& 2k_2u_4+k_3u_1u_8+k_{10}u_9-k_4u_1u_7,  \nonumber\\
  R_8(u_1,...,u_{10}) &=& 4 k_1u_3u_9-k_3u_1u_8, \nonumber\\
  R_9(u_1,...,u_{10}) &=& k_4u_1u_7+2k_8u_{10}-(k_1u_3-k_9u_2+k_{10})u_9,  \nonumber\\
  R_{10}(u_1,...,u_{10}) &=& k_7u_5-k_8u_{10}. \nonumber
\end{eqnarray}
\begin{table}[!htb]\caption {The chemical reactions of the model \label{Tabl1n}}
\begin{center}
\begin{tabular}{|c|l||c|l|}
\hline $1$    & $HC+OH \rightarrow 4RO_2 + 2ALD$      &$6$        & $NO+O_3 \rightarrow NO_2 + O_2$                \\
\hline $2$    & $ALD + h\nu \rightarrow 2HO_2+CO$        &$7$        & $O_3+h \nu \rightarrow O_2 + O(^{1}D)$        \\
\hline $3$    & $RO_2+NO \rightarrow NO_2+ALD+HO_2$     &$8$        & $O(^{1}D) + H_2 O \rightarrow 2 OH$     \\
\hline $4$    & $NO+HO_2 \rightarrow NO_2+OH$            &$9$        & $NO_2+OH \rightarrow HNO_3$                    \\
\hline $5$    & $NO_2+h \nu \rightarrow NO+ O_33$     &$10$       & $CO+OH \rightarrow CO_2+HO_2$                  \\
\hline
\end{tabular}
\end{center}
\end{table}
The chemical part of the model is given in Table \ref{Tabl1n} for the sake of completeness.
The rate coefficients can be found in Table  \ref{Tabl2}.
Some of the coefficients belong to photochemical reactions (the ones with term $h \nu$),
which means that this reactions depend on the light, more precisely on the position of the Sun relative to the horizon:
in $k_2$, $k_5$ and $k_7$ the angle $\theta$ denotes the solar zenith angle,
which is the angle of the Sun measured from vertical.
The chemical species involved in the simplified reactions are written in Table \ref{Tabl3}.
\begin{table}[tb]
\caption {The coefficients of the chemical reactions\label{Tabl2}}
\begin{center}
\begin{tabular}{|c|l||c|l|}
\hline $k_1$    &$6.0e-12$                             &$k_6$         & $1.6e-14$               \\
\hline $k_2$    & $7.8e-05.\exp(-0.87/\cos \theta)$       &$k_7$        & $1.6e-04.\exp(-1.9/\cos \theta)$                      \\
\hline $k_3$    &$8.0e-12$                               &$k_8$        & $2.3e-10$                  \\
\hline $k_4$    & $8.0e-12$                               &$k_9$        & $1.0e-11$                     \\
\hline $k_5$     &$1.0e-02.\exp(-0.39/\cos \theta)$      &$k_{10}$  & $2.9e-13$                    \\
\hline
\end{tabular}
\end{center}
\end{table}
\begin{table}[!htb]\caption {The chemical species in the model\label{Tabl3}}
\begin{center}
\begin{tabular}{|c|c|c|c|c|c|c|c|c|c|}
\hline $u_1$    &$u_2$      &$u_3$         & $u_4$ & $u_5$ & $u_6$ & $u_7$ & $u_8$ & $u_9$ & $u_{10}$               \\
\hline $NO$    & $NO_2$      &$HC$        & $ALD$     & $O_3$      &$HNO_3$        & $HO_2$     & $RO_2$     & $OH$     & $O(^{1} D)$  \\
\hline
\end{tabular}
\end{center}
\end{table}

From both the practical and mathematical point of view, one is naturally interested in the existence and qualitative of the solutions to the problem (\ref{2Dsystem1})-(\ref{starn}). The well-posedness of initial boundary value problems for a system more general than (\ref{2Dsystem1}) is obtained in \cite{Pao}. Throughout of the rest of the paper we assume existence and uniqueness of classical solution of problem  (\ref{2Dsystem1})-(\ref{starn}) which means a function that belongs to
$C([0,T]\times\overline{\Omega}) \bigcap C^1 ((0,T);C(\overline{\Omega})) \bigcap (C(0,T);C^2(\overline{\Omega})) $ and satisfies the equations (\ref{2Dsystem1})-(\ref{2Dsystem3}) pointwise. Moreover, at the finite difference approximations in Sections 3, 4 we assume fourth in time and sixth in space derivatives.

Since we are interested in systems describing chemical concentrations, the nonnegativity of the solutions has to be preserved. It is proved in \cite{CheLi}, that if:
\begin{itemize}
\item[1.] $\mathbf{u}_0(x,y) \geq 0$;
\item[2.] $R_l (x,y,\mathbf{u})$, $l=1,...,L$ is Lipshitz continuous with respect to the concentrations $u_1,u_2,\dots,u_L$ and it satisfies the inequality
   $R_l (x,y,\mathbf{u}) \geq 0$, whenever $u_l=0$, and $\mathbf{u} \in R_{+}^{L} \equiv \{u_k \geq 0,\, k=1,...,L  \},$
\end{itemize}
than $\mathbf{u} \geq 0$ for all $(x,y) \in \Omega$ and $t \in [0,T].$

It is easily to check that the chemical reactions $R_l(u_1,u_2,\dots, u_{10})$, $l=1,...,10$ given by (\ref{starn}) satisfy the point 2. and the solution of problem (\ref{2Dsystem1})-(\ref{2Dsystem3}) with (\ref{starn}) is nonnegative in time $t>0$ if the initial data $\mathbf{u}_0(x,y) \geq 0.$

\section{Central Difference Schemes and Richardson Extrapolation}
In this section, for clarity exposition  we describe the construction of the
second order CDS for the weakly coupled system of two equations
\begin{subequations}\label{2Dsystem2n}
\begin{alignat}{1}
 \frac{\partial u}{\partial t} - a(x,y) \frac{  \partial^2 u  }{   \partial x^2 }
 - b(x,y) \frac{  \partial^2 u  }{   \partial y^2 }
  + c(x,y) \frac{ \partial u }{ \partial x }
+ d(x,y) \frac{ \partial u }{ \partial y }
  &= r (x,y,t,u,v),
  \label{2Dsystem21} \\
 \frac{\partial v}{\partial t} -e(x,y) \frac{  \partial^2 v  }{   \partial x^2 }
 - f(x,y) \frac{  \partial^2 v  }{   \partial y^2 }
  + g(x,y) \frac{ \partial v }{ \partial x }
+ h(x,y) \frac{ \partial v }{ \partial y }
  &= s (x,y,t,u,v),\label{2Dsystem22}
\end{alignat}
\end{subequations}
defined on the cylindric domain $Q_T=\Omega\times(0,T]$, where $\Omega\subset
R^2$ is a bounded domain with Lipshitz boundary. The nonlinear
functions $r$ and $s$ are sufficiently smooth of their arguments.
The coefficients $a(x,y)$,  $b(x,y)$,  $e(x,y)$ and $f(x,y)$ are positive in $\Omega$. We
consider Dirichlet boundary conditions
\begin{equation}\label{BC}
    u(x,y,t) = \bar{\phi}(x,y,t), \;\;v(x,y,t) = \bar{\bar{\phi}}(x,y,t),\;\; (x,y,t)\in
    \partial\Omega\times(0,T]
\end{equation}
and initial conditions
\begin{equation}\label{IC}
    u(x,y,0) = \bar{\psi}(x,y), \;\;v(x,y,0) = \bar{\bar{\psi}}(x,y),\;\;
    (x,y)\in\Omega ,
\end{equation}
where $\bar{\phi}$, $\bar{\bar{\phi}}$, $\bar{\psi}$ and
$\bar{\bar{\psi}}$ are given and smooth data and compatibility of the boundary and initial data is ensured.

Let for simplicity
the domain $\Omega $ is a rectangle $\Omega=[0,X]\times[0, Y]$.  We
introduce uniform meshes in the following way:  $ \overline{\omega}_{h,x} = \{ x_{i} =ih_x
, \; \; i=0,1, \dots , M_x, \; \; h_x = X/M_x \} $, $
\overline{\omega}_{h,y}  = \{ y_{j} =jh_y , \; \; j=0,1, \dots ,M_y,
\; \; h_y = Y/M_y \} $ and then
$\overline{\Omega}_h=\omega_{h,x}\times\omega_{h,y}$,
$\overline{\Omega}_h={\Omega}_h\cup\partial\Omega_h$, where
$\Omega_h$ consist of all interior mesh points and
$\partial\Omega_h$ - of all boundary mesh points.

We will used the index pair $(i,j)$ to represent the mesh point
$(x_i,y_j)$ and define
\begin{equation*}
    u_{i,j}=u(x_i,y_j,t),
    \;\;v_{i,j}=v(x_i,y_j,t),\;\;r_{i,j}=r(x_i,y_j,t,u_{i,j},v_{i,j}),
    \;\;\; \mbox{ect.}
\end{equation*}
For $w=u, v$ we introduce the central difference operators
\begin{equation}\label{DO}
\begin{array}{lcl}
  \delta_x w_{i,j}=(w_{i+1,j}-w_{i-1,j})/(2h_x) , &&
 \delta^2_x w_{i,j}=(w_{i+1,j}-2w_{i,j}+w_{i-1,j})/h^2_x,\\ \delta_y w_{i,j}=(w_{i+1,j}-w_{i-1,j})/(2h_y) ,
  && \delta^2_y
  w_{i,j}=(w_{i,j+1}-2w_{i,j}+w_{i,j-1})/h^2_y.
\end{array}
\end{equation}
\subsection{Second-order space semidiscretization}
Application of the difference operators (\ref{DO}) into the
system (\ref{2Dsystem2n})  for every point $(i,j)\in\Omega_h$ leads to
\begin{eqnarray*}\label{SD}
    \left.\frac{\partial u}{\partial t}\right|_{(x_i,y_j)}-a_{i,j}\delta^2_x u_{i,j}-b_{i,j}\delta^2_y u_{i,j}
    +c_{i,j}\delta_x u_{i,j}+d_{i,j}\delta_y u_{i,j}+\chi_{i,j,1}&=&
    r_{i,j},\\
    \left.\frac{\partial v}{\partial t}\right|_{(x_i,y_j)} -e_{i,j}\delta^2_x v_{i,j}-f_{i,j}\delta^2_y v_{i,j}
    +g_{i,j}\delta_x v_{i,j}+h_{i,j}\delta_y v_{i,j}+\chi_{i,j,2}&=&
   s_{i,j},
\end{eqnarray*}
where the truncation errors $\chi_{i,j,1}$ and  $\chi_{i,j,2}$ are
\begin{equation}\label{TE}
\begin{array}{lcl}
    \chi_{i,j,1}&=&\frac{h_x^2}{12}\left(2c\frac{\partial^3u}{\partial x^3}
    -a \frac{\partial^4u}{\partial x^4}\right)_{i,j}+
    \frac{h_y^2}{12}\left(2d\frac{\partial^3u}{\partial y^3}
    -b \frac{\partial^4u}{\partial
    y^4}\right)_{i,j}+\mathcal{O}(h^4_x+h^4_y), \\
    \chi_{i,j,2}&=&\frac{h_x^2}{12}\left(2g\frac{\partial^3v}{\partial x^3}
    -e \frac{\partial^4v}{\partial x^4}\right)_{i,j}+
    \frac{h_y^2}{12}\left(2h\frac{\partial^3v}{\partial y^3}
    -f \frac{\partial^4v}{\partial
    y^4}\right)_{i,j}+\mathcal{O}(h^4_x+h^4_y).
\end{array}
\end{equation}
After dropping the truncation error terms a semi-discrete second-order central difference approximation of (\ref{2Dsystem2n}) is obtained:
\begin{equation}\label{SDN}
\begin{array}{lcl}
    \left.\frac{\partial u^h}{\partial t}\right|_{(x_i,y_j)}-a_{i,j}\delta^2_x u_{i,j}^h-b_{i,j}\delta^2_y u_{i,j}^h
    +c_{i,j}\delta_x u_{i,j}^h+d_{i,j}\delta_y u_{i,j}^h&=&
    r_{i,j}^h,\\
    \left.\frac{\partial v^h}{\partial t}\right|_{(x_i,y_j)} -e_{i,j}\delta^2_x v_{i,j}^h-f_{i,j}\delta^2_y v_{i,j}^h
    +g_{i,j}\delta_x v_{i,j}^h+h_{i,j}\delta_y v_{i,j}^h&=&
   s_{i,j}^h,
\end{array}
\end{equation}
where for $ (i,j)\in\Omega_h$
 \begin{eqnarray*}
u_{ i,j}^h\approx u(x_i,y_j,t), \quad&& v_{i,j}^h\approx v(x_i,y_j,t), \\
r_{i,j}^h\approx r(x_i,y_j,t,u_{i,j}^h,v_{i,j}^h),\quad && s_{i,j}^h\approx s(x_i,y_j,t,u_{i,j}^h,v_{i,j}^h).
\end{eqnarray*}

Now we introduce the matrix representation for the system
(\ref{SDN}). We order the mesh points
lexicographically from left to right in $x$ direction and from
the bottom to the top in $y$ direction. Excluding the boundary mesh
points $(i,j)\in\partial\Omega_h$, for $j=1,2,...,M_y-1$ we define
the following $(M_x-1)$ dimensional vectors:
\begin{eqnarray*}
  U_{j}^h=\left(u_{1,j}^h, u_{2,j}^h,...,u_{M_x-1,j}^h
 \right),\;\;\;&&
 V_{j}^h=\left(v_{1,j}^h, v_{2,j}^h,...,v_{M_x-1,j}^h  \right), \\
 R_j(U_{j}^h,V_{j}^h)=\left(
 R_{1,j},R_{2,j},...,R_{M_x-1,j}\right),&&
 S_j(U_{j}^h,V_{j}^h)=\left(
 S_{1,j},\;S_{2,j},...,S_{M_x-1,j}\right)\;\;\;
\end{eqnarray*}
and then
\begin{eqnarray*}
   U=\left(U_{1}^h, U_{2}^h,...,U_{M_y-1}^h
 \right)^T,\;\;\;&& V=\left(V_{1}^h, V_{2}^h,...,V_{M_y-1}^h
 \right)^T, \\
  R=\left(R_{1}, R_{2},...,R_{M_y-1}\right)^T,\;\;\;
 &&S=\left(S_{1}, S_{2},...,S_{M_y-1}\right)^T.\;\;\;
\end{eqnarray*}

We then rewrite the system (\ref{SDN}) as a system of ordinary differential equations
\begin{eqnarray}\label{systemPQuvn}
 \frac{d}{dt}U+\bar{P}U&=&R+\bar{\Phi}, \quad t\in(0,T], \label{systemPQun}\\
  \frac{d}{dt}V+\bar{\bar{P}}V&=&S+\bar{\bar{\Phi}},\quad t\in(0,T] \label{systemPQvn}
\end{eqnarray}
with initial conditions $U(0)$ and $V(0)$ obtaining from
$\bar{\psi}$ and $\bar{\bar{\psi}}$ for $(i,j)\in\Omega_h$ after the
reordering. In (\ref{systemPQun}) the
matrix $\bar{P}$ is
$(M_y-1)\times(M_y-1)$ block-tridiagonal matrix
$\bar{P}=tridiag(\bar{P}_{k,k-1},\bar{P}_{k,k},\bar{P}_{k,k+1})$ and
$\bar{P}_{k,l}$, $l=k-1,k,k+1$ are tridiagonal matrixes for $l=k$ and diagonal for $l=k\pm1$ of
order $(M_x-1)\times(M_x-1)$. Let for two natural numbers  $m$ and
$M$, $m<M$ denote $m:M=m,m+1,...,M$ and assume that $\mathbf{p}_{k,m:M}$ is a vector with entrances $\mathbf{p}_{k,m:M}=(p_{k,m},p_{k,m+1},...,p_{k,M})$ . Then from (\ref{SDN}) and (\ref{DO}) the entrances
of $\bar{P}_{k,l}$ are
\begin{equation}\label{barPstand}
\bar{P}_{k,l}=tridiag(\mathbf{p}_{k,2:M_x-1}^{(-1,\varepsilon)},
\mathbf{p}_{k,2:M_x}^{(0,\varepsilon)},\mathbf{p}_{k,1:M_x-2}^{(1,\varepsilon)})\;\;\;\quad
l=k+\varepsilon, \;\; \varepsilon=0,\pm1 \, ,
\end{equation}
where
\begin{eqnarray}\label{Pstand}
\nonumber    p_{i,j}^{(\pm1,0)}&=& \pm\frac{c(i,j)}{2h_x}-\frac{a(i,j)}{h_x^2} \,, \\
                      p_{i,j}^{(0,\pm1)}&=& \pm\frac{d(i,j)}{2h_y}-\frac{b(i,j)}{h_y^2} \,,\\
\nonumber    p_{i,j}^{(0,0)}&=&2\frac{a(i,j)}{h_x^2}+2\frac{b(i,j)}{h_y^2}\,.
\end{eqnarray}
Replacing $a\leftrightarrow e$, $b\leftrightarrow f$, $c\leftrightarrow g$ and $d\leftrightarrow h$ in a similar way we obtain the entrances of the matrix $\bar{\bar{P}}$.

The vectors $\bar{\Phi}$ and $\bar{\bar{\Phi}}$ in (\ref{systemPQun})-(\ref{systemPQvn}) are associated with
the boundary functions and also depend on time $t$.
\subsection{Full discretization}
For discretization in time the so called $\theta$-weight method is used. Let $
\omega_{\tau }  = \{ t_{n} =n \tau  , \; \; n=0,1, \dots , N, \; \;
\tau = T/N \} $ be uniform mesh in time with time step $\tau$. Then the
weight $\theta$-discretization of (\ref{systemPQun}), (\ref{systemPQvn}) may be written in the following way:,
\begin{eqnarray}\label{fulldiscr}
 \frac{U^{n+1}-U^{n}}{\tau}+\bar{P}U^{n,\theta}&=&R^{n,\theta}+\bar{\Phi}^{n,\theta}, \quad t\in(0,T), \\
  \frac{V^{n+1}-V^{n}}{\tau}+\bar{\bar{P}}V^{n,\theta}&=&S^{n,\theta}+\bar{\bar{\Phi}}^{n,\theta},\quad t\in(0,T), \nonumber
\end{eqnarray}
where
$Z^{n,\theta}=\theta Z^{n+1}+(1-\theta )Z^{n}$ for $Z=U,V,R,S,\bar{\Phi},\bar{\bar{\Phi}}$, $Z^n\approx Z(t_n)$ and $0\leq\theta\leq 1$, $n=0,1, \dots N-1$. For $\theta=1$ one obtain the fully implicit finite difference scheme,
for $\theta=0$ - explicit and for $\theta=1/2$ - the Crank-Nicolson scheme. The last case has an advantage that the scheme is of second order in time and as we  want to derive schemes of higher order, in the numerical experiments we use mainly $\theta=1/2$.

For $\theta >0$ the finite difference schemes requires solving of
nonlinear algebraic systems. We briefly discuss the application of
the Newton method on the problem (\ref{fulldiscr}). To apply the
classical Newton method the system (\ref{fulldiscr}) is rewritten in
the form $\Upsilon ( { W } ) =0$,  where $ {W } = [
U ^T, V^T ]^T $ is a vector of length $2(M_x-1)(M_y-1)$.
We set $\stackrel{0\;\;\;}{W^{n+1}}$ as initial guess on the new time layer $t=t_{n+1}$ to be the numerical solution on the
previous time layer $t=t_{n}$. Then to find the solution on
$t=t_{n+1} $ the iterative process with appropriate stopping
criteria is used:
\begin{equation} \label{Newton1}
   \left\{
   \begin{array}{l}
   \Upsilon ' (\stackrel{k\;\;\;\;}{W^{n+1}} )  \stackrel{k}{\Delta} = - \Upsilon (\stackrel{k\;\;\;\;}{W^{n+1}} )\,,  \\
   \stackrel{k+1\;\;\;\;}{W^{n+1}} =\stackrel{k\;\;\;\;}{W^{n+1}} + \stackrel{k}{\Delta}\,.
   \end{array}
   \right.
   \end{equation}
Here $ \stackrel{k}{\Delta}$ is a vector of the increments and the Jacobian matrix $\Upsilon'  (\stackrel{k\;\;\;}{W^{n+1}} )$ for $\theta=1/2$ is
\begin{equation}
\Upsilon'  (\stackrel{k\;\;\;}{W^{n+1}} )=\frac{\partial\Upsilon}{\partial {W}}=
 \left.
 \left(
 \begin{array}{ c   c }
 \frac{1}{\tau}I + \frac{1}{2}\bar{P}- \frac{1}{2} \frac{  \partial R }{ \partial U }&  \;\; \frac{1}{2} \frac{  \partial R }{   \partial V  }\\
  \frac{1}{2}\frac{\partial S }{ \partial U }& \;\;  \frac{1}{\tau}I + \frac{1}{2}\bar{\bar{P}}- \frac{1}{2} \frac{  \partial S }{ \partial V }
 \end{array}
 \right)
 \right\vert_{  (U,V) = (\stackrel{k}{U},\stackrel{k}{V} )  } ,
 \end{equation}
where $I$ is the identity matrix and $\overline{P}$, $\overline{\overline{P}}$ - as defined by (\ref{barPstand}), (\ref{Pstand}).
In the numerical experiments to solve the first line in (\ref{Newton1})  which is a linear system of $2(M_x-1)(M_y-1)$ equations we use the so called inexact Newton method \cite{DemSta}, i.e. we solve this system approximately using the MatLab function bicgstab(l) (biconjugate gradients stabilized (l) method) that gives better results for our examples in sense of convergence of the inner iterations and the CPU time.
\subsection{Richardson extrapolation}
Richardson extrapolation is a powerful  computational tool which can successfully be used in the efforts to improve the accuracy of the of the approximate solutions of the systems of partial differential equations (PDEs) obtained by finite difference methods.

Therefor, another way for obtaining the difference schemes of higher order is to use the Richardson extrapolation method.
The main idea \cite{MarSha}  is to solve the difference scheme on two or more consecutive meshes and then to combine the obtained numerical solutions with appropriate weights.
Let assume that $h_x=h_y=h$ and for the numerical solution on the $n$-th time layer the following expression is true:
\begin{equation}\label{RE}
    U^{\tau}_h=U^{n}_{(i,j)}=u(x_i,y_j,t^n)+C_1h^\sigma+\chi(h,\tau),\quad  (x_i,y_j,t_n) \in \Omega_{h,\tau},
\end{equation}
where function $\chi(h,\tau)$ is a remainder term and $C_1$ does not depend on $h_x$, $h_y$ and $\tau$. If we want to eliminate the term $C_1h^\sigma$, we do the following steps:
\begin{itemize}
  \item solve the difference scheme on two consecutive meshes: coarse  one $\Omega_{h,\tau }$ and fine one $\Omega_{h/2,\tau }$ and let the corresponding numerical solutions be $U^{\tau}_h $ and   $U^{\tau}_{h/2} $;
  \item  find the weights $ \gamma_1$ and $\gamma_2$ from the system
  \begin{eqnarray}\label{REsys}
    \gamma_1+\gamma_2&=& 1 \\
     \gamma_1+\frac{\gamma_2}{2^\sigma}&=& 0 \nonumber
  \end{eqnarray}
    \item obtain a new numerical solution on the coarse mesh
     \begin{equation*}
  U_{extr}=\gamma_1U^{\tau}_h+\gamma_2U^{\tau}_{h/2} \qquad (x_i,y_j,t_n)\in\Omega_{h,\tau }\; .
 \end{equation*}
 \end{itemize}
From (\ref{REsys}) we have for the case of central Crank-Nicolson Scheme ($\sigma=2$) that the coefficients for the Richardson extrapolation are
\begin{equation}\label{RE2}
   \gamma_1=-1/3 \qquad\gamma_2=4/3.
\end{equation}
In the case of CFDS and Richardson Extrapolation ($\sigma=4$) the corresponding weight coefficients are
\begin{equation}\label{RE4}
   \gamma_1=-1/15 \qquad\gamma_2=16/15.
\end{equation}
If in (\ref{RE}) the more detailed analysis of the LTE is done, then the prolongation of the idea of space-time Richardson extrapolation \cite{Rich} can be applied.

\section{Compact Difference Schemes}
In this section, just for clarity  we describe the construction of the
CFDS  again for the system of two equations (\ref{2Dsystem2}).
\subsection{Space discretization}

In order to eliminate the terms of $\mathcal{O}(h^2_x+h^2_y)$ in
(\ref{TE}) we differentiate the equation (\ref{2Dsystem21}) twice
with respect to $x$ obtaining  expressions for
$\frac{\partial^3u}{\partial x^3}$, $\frac{\partial^4u}{\partial
x^4}$, and twice with respect to $y$ for
$\frac{\partial^3u}{\partial y^3}$, $\frac{\partial^4u}{\partial
y^4}$.

Let
\begin{equation*}
    \tilde{a}_{i,j}=(c_{i,j}+2\delta_x a_{i,j})/a_{i,j}, \;\;\;\;
    \tilde{b}_{i,j}=(d_{i,j}+2\delta_y b_{i,j})/b_{i,j}, \;\;\;\;
    (i,j)\in\Omega_h\;.
\end{equation*}
Let also
\begin{eqnarray*}
  \alpha_{i,j} &=& a_{i,j}+\frac{h_x^2}{12}
  \left(\delta_x^2a_{i,j}- \tilde{a}_{i,j}(\delta_xa_{i,j} -c_{i,j})-2\delta_x c_{i,j}   \right)
  + \frac{h_y^2}{12}\left(\delta_y^2a_{i,j}- \tilde{b}_{i,j}\delta_ya_{i,j}\right) ,\\
  \beta_{i,j} &=& b_{i,j}+\frac{h_x^2}{12}
  \left(\delta_x^2b_{i,j}- \tilde{a}_{i,j}\delta_xb_{i,j}\right)
  +\frac{h_y^2}{12}\left(\delta_y^2b_{i,j}-
  \tilde{b}_{i,j}(\delta_yb_{i,j}-d_{i,j})-2\delta_yd_{i,j}\right) , \\
  \tilde{\alpha}_{i,j} &=&c_{i,j}+\frac{h_x^2}{12}
  \left(\delta_x^2c_{i,j}- \tilde{a}_{i,j}\delta_xc_{i,j}\right)
  +\frac{h_y^2}{12}\left(\delta_y^2c_{i,j}-
  \tilde{b}_{i,j}\delta_yc_{i,j}\right) ,  \\
  \tilde{\beta}_{i,j} &=&d_{i,j}+\frac{h_x^2}{12}
  \left(\delta_x^2d_{i,j}- \tilde{a}_{i,j}\delta_xd_{i,j}\right)
  +\frac{h_y^2}{12}\left(\delta_y^2d_{i,j}-
  \tilde{b}_{i,j}\delta_yd_{i,j}\right) ,
\end{eqnarray*}
and
\begin{eqnarray*}
  \theta_{i,j} &=& \frac{h_y^2}{12}c_{i,j}
   -\frac{h_x^2}{12}(2\delta_xb_{i,j}-\tilde{a}_{i,j}b_{i,j}), \;\;\;
   \tilde{\theta}_{i,j} = \frac{h_x^2}{12}d_{i,j}
   -\frac{h_y^2}{12}(2\delta_ya_{i,j}-\tilde{b}_{i,j}a_{i,j}),  \\
  \gamma_{i,j}
  &=&\frac{h_x^2}{12}b_{i,j}+\frac{h_y^2}{12}a_{i,j},\;\;\;
  \tilde{\gamma}_{i,j}=\frac{h_x^2}{12}(2\delta_x-\tilde{a}_{i,j}d_{i,j})
  +\frac{h_y^2}{12}(2\delta_yc_{i,j}-\tilde{b}_{i,j}c_{i,j}).
\end{eqnarray*}
Define the following difference operators
\begin{eqnarray*}
 l_{i,j}^h&=&-\alpha_{i,j}\delta_x^2-\beta_{i,j}\delta_y^2 +\tilde{\alpha}_{i,j}\delta_x
 +\tilde{\beta}_{i,j}\delta_y-\gamma_{i,j}\delta_x^2 \delta_y^2
 +\theta_{i,j}\delta_x\delta_y^2+\tilde{\theta}_{i,j}\delta_x^2 \delta_y
 +\tilde{\gamma}_{i,j}\delta_x \delta_y\\
 \nu_{i,j}^h&=&1+\frac{h_x^2}{12}(\delta_x^2-\tilde{a}_{i,j}\delta_x)
 +\frac{h_y^2}{12}(\delta_y^2-\tilde{b}_{i,j}\delta_y).
 \end{eqnarray*}
Applying these operators to (\ref{2Dsystem21}) we have
\begin{equation}\label{operator}
l_{i,j}^hu_{i,j}=\nu_{i,j}^h(r_{i,j}-u_{t,i,j})+\mathcal{O}(h_x^4+h_x^2h_y^2+h^4_y).
\end{equation}
For convenience,  we introduce also the
operators
\begin{equation*}\label{Operator}
    {\bar{\mathcal{P}}}_{i,j}^h=6h_x^2l_{i,j}^h, \;\;\; {\bar{\mathcal{Q}}}_{i,j}^h=6h_x^2\nu_{i,j}^h.
\end{equation*}
Let $\sigma=h_x/h_y$ be the ratio of the mesh sizes. Then
\begin{eqnarray*}\label{OPERATOR}
 {\bar{\mathcal{P}}}_{i,j}^hu_{i,j}&=&\sum_{k_1=-1}^1\sum_{k_2=-1}^1p_{i,j}^{(k_1,k_2)}u_{i+k_1,j+k_2},\;\;\;\\
{\bar{\mathcal{Q}}}_{i,j}^hu_{i,j}&=&\sum_{k_1=-1}^1\sum_{k_2=-1}^1q_{i,j}^{(k_1,k_2)}u_{i+k_1,j+k_2},
\end{eqnarray*}
where
\begin{eqnarray}\label{P}
 \nonumber p_{i,j}^{(\pm1,-1)} &=&-\frac{{a}_{i,j}+\sigma^2{b}_{i,j}}{2}
 \pm\frac{1}{4} \left(c_{i,j}-\sigma^2(2\delta_xb_{i,j}-\tilde{a}_{i,j}b_{i,j})
 \mp\sigma d_{i,j}\pm\frac{1}{\sigma} (2\delta_ya_{i,j}-\tilde{b}_{i,j}a_{i,j}) \right)h_x   \\
 \nonumber &&\mp\frac{1}{8}\left(\sigma (2\delta_x-\tilde{a}_{i,j}d_{i,j})
 +\frac{1}{\sigma} (2\delta_yc_{i,j}-\tilde{b}_{i,j}c_{i,j})  \right)h_x^2,  \\
 \nonumber p_{i,j}^{(\pm1,1)} &=&-\frac{{a}_{i,j}+\sigma^2{b}_{i,j}}{2}
 \pm\frac{1}{4} \left(c_{i,j}-\sigma^2(2\delta_xb_{i,j}-\tilde{a}_{i,j}b_{i,j})
 \pm\sigma d_{i,j}\mp\frac{1}{\sigma} (2\delta_ya_{i,j}-\tilde{b}_{i,j}a_{i,j}) \right)h_x   \\
 \nonumber &&\pm\frac{1}{8}\left(\sigma (2\delta_x-\tilde{a}_{i,j}d_{i,j})
 +\frac{1}{\sigma} (2\delta_yc_{i,j}-\tilde{b}_{i,j}c_{i,j})  \right)h_x^2,  \\
  p_{i,j}^{(\pm1,0)}&=&\sigma^2b_{i,j}-5a_{i,j}\pm\left(3\tilde{\alpha}_{i,j}-\frac{1}{2}c_{i,j}+
  \frac{\sigma^2}{2}(2\delta_yc_{i,j}-\tilde{b}_{i,j}c_{i,j})     \right)h_x  \\
 \nonumber &&-\frac{1}{2}\left(\delta_x^2a_{i,j}- \tilde{a}_{i,j} (\delta_xa_{i,j}-c_{i,j})
 -2\delta_xc_{i,j}+\frac{1}{\sigma^2}(\delta^2_ya_{i,j}-\tilde{b}_{i,j}\delta_ya_{i,j})  \right)h_x^2  \\
 \nonumber p_{i,j}^{(0,\pm1)}&=&a_{i,j}-5\sigma^2b_{i,j}
 \pm\left(3\sigma\tilde{\beta}_{i,j}-\frac{\sigma}{2}d_{i,j}
 +\frac{1}{2\sigma}(2\delta_ya_{i,j}-\tilde{b}_{i,j}a_{i,j})\right)h_x\\
 \nonumber&&-\frac{1}{2}\left(\sigma^2(\delta_x^2b_{i,j}-\tilde{a}_{i,j}\delta_xb_{i,j})
  +\delta_y^2b_{i,j}-2\delta_yd_{i,j}-\tilde{b}_{i,j}(\delta_yb_{i,j}-d_{i,j})
  \right)h_x^2,\\
\nonumber p_{i,j}^{(0,0)}&=&10(a_{i,j}+\sigma^2b_{i,j})+ \left(
 \delta_x^2a_{i,j}-\tilde{a}_{i,j}(\delta_xa_{i,j}-c_{i,j})-2\delta_xc_{i,j}
 +\frac{1}{\sigma^2}(\delta^2_ya_{i,j}-\tilde{b}_{i,j}\delta_ya_{i,j})
 \right)h_x^2\\
\nonumber&&+\left(\sigma^2(\delta_x^2b_{i,j}-\tilde{a}_{i,j}\delta_xb_{i,j})
+\delta_y^2b_{i,j}-2\delta_yd_{i,j}-\tilde{b}_{i,j}(\delta_yb_{i,j}-d_{i,j})
\right)h_x^2
\end{eqnarray}
and
\begin{equation}\label{Q}
q_{i,j}^{(\pm1,\pm1)}=0,\;q_{i,j}^{(\pm1,0)}=\frac{1}{4}(2\mp\tilde{a}_{i,j}h_x)h_x^2,
\;q_{i,j}^{(0,\pm1)}=\frac{1}{4}(2\mp\frac{\tilde{a}_{i,j}}{\sigma}h_x)h_x^2,
\;q_{i,j}^{(0,0)}=4h_x^2 .
\end{equation}
With these notations, after dropping the term
$\mathcal{O}(h_x^4+h_x^2h_y^2+h^4_y)$ in (\ref{operator}) the
semi-discrete compact finite difference approximation of (\ref{2Dsystem21}) and the left parts of (\ref{BC}),  (\ref{IC}) are as
follows:
\begin{equation}\label{SDFDu}
\left\{\begin{array}{lcc}
  {\bar{\mathcal{P}}}_{i,j}^h u_{i,j}^h={\bar{\mathcal{Q}}}_{i,j}^h\left(
 r_{i,j}^h-\frac{d}{dt}u_{i,j}^h\right),& \;\;(i,j)\in\Omega_h,\;\;& t\in(0,T],\\
\quad\;\;\;u_{i,j}^h=\bar{\phi}_{i,j},\;\;\;&\;\;\;(i,j)\in\partial\Omega_h ,\;\;&t\in(0,T],\\
\quad\;\;\;u_{i,j}^h=\bar{\psi}_{i,j},\;\;\;&\;\;\;(i,j)\in\bar\Omega_h,\;\;&t=0.
\end{array}\right.
\end{equation}
In a similar way we proceed with the equation (\ref{2Dsystem22}) and the left parts of (\ref{BC}),  (\ref{IC}).
Replacing $a_{i,j}$, $b_{i,j}$, $c_{i,j}$, $d_{i,j}$
with $e_{i,j}$, $f_{i,j}$, $g_{i,j}$, $h_{i,j}$
and $\bar{\mathcal{P}}_{i,j}^h $, $\bar{\mathcal{Q}}_{i,j}^h $ with
$\bar{\bar{\mathcal{P}}}_{i,j}^h $, $\bar{\bar{\mathcal{Q}}}_{i,j}^h
$ we obtain the second part of the semi-discrete  nonlinear system
\begin{equation}\label{SDFDv}
\left\{\begin{array}{lcc}
 \bar{\bar{\mathcal{P}}}_{i,j}^h v_{i,j}^h=\bar{\bar{\mathcal{Q}}}_{i,j}^h\left(
 s_{i,j}^h-\frac{d}{dt}v_{i,j}^h\right),& \;\;(i,j)\in\Omega_h,\;\;& t\in(0,T],\\
\quad\;\;\; v_{i,j}^h=\bar{\bar{\phi}}_{i,j},\;\;\;
&\;\;\;(i,j)\in\partial\Omega_h ,\;\;&t\in(0,T],\\
\quad\;\;\;
v_{i,j}^h=\bar{\bar{\psi}}_{i,j},\;\;\;&\;\;\;(i,j)\in\bar\Omega_h,\;\;&t=0.
\end{array}\right.
\end{equation}

Now we introduce the matrix representation for the system
(\ref{SDFDu}), (\ref{SDFDv}).
We obtain the following system of ordinary differential equations
\begin{eqnarray}
  \bar{Q}\frac{d}{dt}U^h+\bar{P}U^h&=&\bar{Q}R+\bar{\Phi}, \quad t\in(0,T], \label{systemPQu}\\
  \bar{\bar{Q}}\frac{d}{dt}V^h+\bar{\bar{P}}V^h&=&\bar{\bar{Q}}S+\bar{\bar{\Phi}} \label{systemPQv}
\end{eqnarray}
with initial conditions $U^h(0)$ and $V^h(0)$ obtaining from
$\bar{\psi}$ and $\bar{\bar{\psi}}$ for $(i,j)\in\Omega_h$ after the
reordering. In system (\ref{systemPQu}), (\ref{systemPQv}) the
matrix $\bar{P}$ (similarly $\bar{\bar{P}}$)  is
$(M_y-1)\times(M_y-1)$ block-tridiagonal matrix
$\bar{P}=tridiag(\bar{P}_{k,k-1},\bar{P}_{k,k},\bar{P}_{k,k+1})$ and
$\bar{P}_{k,l}$, $l=k-1,k,k+1$ are also tridiagonal matrixes of
order $(M_x-1)\times(M_x-1)$.  Then from (\ref{P}) the entries
of $\bar{P}_{k,l}$ are
\begin{equation}\label{barP}
\bar{P}_{k,l}=tridiag(p_{k,2:M_x-1}^{(-1,\varepsilon)},
p_{k,2:M_x}^{(0,\varepsilon)},p_{k,1:M_x-2}^{(1,\varepsilon)})\;\;\;\quad
l=k+\varepsilon, \;\; \varepsilon=0,\pm1 \;.
\end{equation}
The entries of $\bar{Q}_{k,l}$ (similarly $\bar{\bar{Q}}$)) are
\begin{equation}\label{barQ}
\bar{Q}_{k,l}=tridiag(q_{k,2:M_x-1}^{(-1,\varepsilon)},\;
q_{k,2:M_x}^{(0,\varepsilon)},\;q_{k,1:M_x-2}^{(1,\varepsilon)})\;\;\;\quad
l=k+\varepsilon, \;\; \varepsilon=0,\pm1
\end{equation}
with a remark that for $\varepsilon=\pm1$ matrixes $\bar{Q}_{k,l}$
are diagonal (instead tridiagonal) matrixes, see (\ref{Q}).

The vectors $\bar{\Phi}$ and $\bar{\bar{\Phi}}$ are associated with
the boundary functions and also depend on time $t$.
\subsection{Time discretization}
For discretization of the ODE system (\ref{systemPQu})-(\ref{barQ}) in time the $\theta$-weight method with $\theta=1/2$ is used in the numerical experiments.  Then the
Crank-Nicolson full discretization of (\ref{systemPQu}), (\ref{systemPQv}) is as follows:
\begin{equation}\label{fulldiscrCFDS}
\begin{array}{lcl}
\bar{Q} \frac{U^{n+1}-U^{n}}{\tau}+\bar{P}U^{n,\theta}&=&\bar{Q}R^{n,\theta}+\bar{\Phi}^{n,\theta}, \quad n=1,...,N-1, \\
 \bar{\bar{Q}} \frac{V^{n+1}-V^{n}}{\tau}+\bar{\bar{P}}V^{n,\theta}&=&\bar{\bar{Q}}S^{n,\theta}+\bar{\bar{\Phi}}^{n,\theta},\quad n=1,...,N-1.
\end{array}
\end{equation}
Similarly to the previous Section we apply the
classical Newton method. The system (\ref{fulldiscrCFDS}) is rewritten in
the form $\Upsilon ( { W } ) =0$,  where $ {W } = [
U ^T, V^T ]^T $ is a vector of length $2(M_x-1)(M_y-1)$.
We set $\stackrel{0\;\;\;}{W^{n+1}}$ as initial guess on the new time layer $t=t_{n+1}$ to be the numerical solution on the
previous time layer $t=t_{n}$. Then to find the solution on
$t=t_{n+1} $ the iterative process with appropriate stopping
criteria is used:
\begin{equation*} \label{Newton1CFDS}
   \left\{
   \begin{array}{l}
   \Upsilon ' (\stackrel{k\;\;\;\;}{W^{n+1}} )  \stackrel{k}{\Delta} = - \Upsilon (\stackrel{k\;\;\;\;}{W^{n+1}} )\,,  \\
   \stackrel{k+1\;\;\;\;}{W^{n+1}} =\stackrel{k\;\;\;\;}{W^{n+1}} + \stackrel{k}{\Delta}\,.
   \end{array}
   \right.
   \end{equation*}
Here $ \stackrel{k}{\Delta}$ is a vector of the increments and the Jacobian matrix $\Upsilon'  (\stackrel{k\;\;\;}{W^{n+1}} )$ for $\theta=1/2$ now is
\begin{equation*}
\Upsilon'  (\stackrel{k\;\;\;}{W^{n+1}} )= \left.
 \left(
 \begin{array}{ c   c }
 \frac{1}{\tau}\bar{Q} + \frac{1}{2}\bar{P}- \frac{1}{2}\bar{Q} \frac{  \partial R }{ \partial U }&  \;\; \frac{1}{2}\bar{Q} \frac{  \partial R }{   \partial V  }\\
  \frac{1}{2}\bar{\bar{Q}}\frac{\partial S }{ \partial U }& \;\;  \frac{1}{\tau}\bar{\bar{Q}} + \frac{1}{2}\bar{\bar{P}}- \frac{1}{2}\bar{\bar{Q}} \frac{  \partial S }{ \partial V }
 \end{array}
 \right)
 \right\vert_{  (U,V) = (\stackrel{k}{U},\stackrel{k}{V} )  } .
 \end{equation*}

\section{Numerical results}
In this section we consider two examples to illustrate the properties of the numerical schemes derived. The first one is an artificial problem with analytical solution and the second one is the two dimensional air-pollution model described in Section 2.
\subsection{Example 1 (known analytical solution)}
Here we consider a problem slightly different  from the problem (\ref{2Dsystem1})-(\ref{starn}):
\begin{equation*}\label{2Dsystem1exact}
 \frac{\partial u_l}{\partial t}
 - K \triangle u_{l} +  \mathbf{b}_{l} . \nabla u_{l}
= R_{l} (x,y, \mathbf{u})+\xi_{l}(x,y,t)  ,  (x,y,t) \in \Omega\times (0,T].
\end{equation*}
The functions  $\xi_l$, $l=1,...,10,$  and the initial and boundary conditions are chosen so that the exact solution is
\begin{equation*}
u_l=\exp(-t/T) sin(\frac{\pi x}{X})sin(\frac{\pi y}{Y}), \quad l=1,...,10, \quad (x,y,t) \in  \overline{\Omega} \times [0,T].
\end{equation*}
The other parameters are as follows: $X=Y=500$, $T=1440$, $\mu=2\pi/(60T)$, $K=1.8$.

For the $l^{th}$ substances with $error_{M,l}$ we denote the
error (the difference between the exact and the numerical solution) in maximum norm, obtained on the last time layer $t_N=T$ for the number of space subintervals $M_x=M_y=M$:
\begin{equation*}
    error_{M,l}=\max_{i,j\in \bar{\Omega}_h}\|u_l(x_i,y_j,t_N)-u_l^h(i,j,N) \|.
\end{equation*}
The ratio between the errors obtained on two consecutive mesh refinements (usually doubling) is denoted by $ratio$:
 \begin{equation*}
    ratio=ratio_{M,l/2M,l}=:error_{M,l}/ error_{2M,l}.
\end{equation*}

\begin{figure}[tbp]
\centering
\psfig{figure=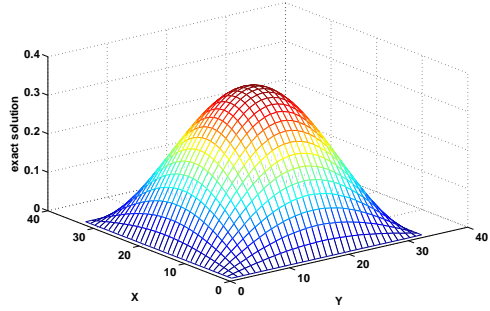,width=7.2cm}
\caption{The exact solution}\label{Fig1}
\end{figure}

\begin{figure}[tbp]
\centering
\begin{tabular}{cc}
\makebox{\psfig{figure=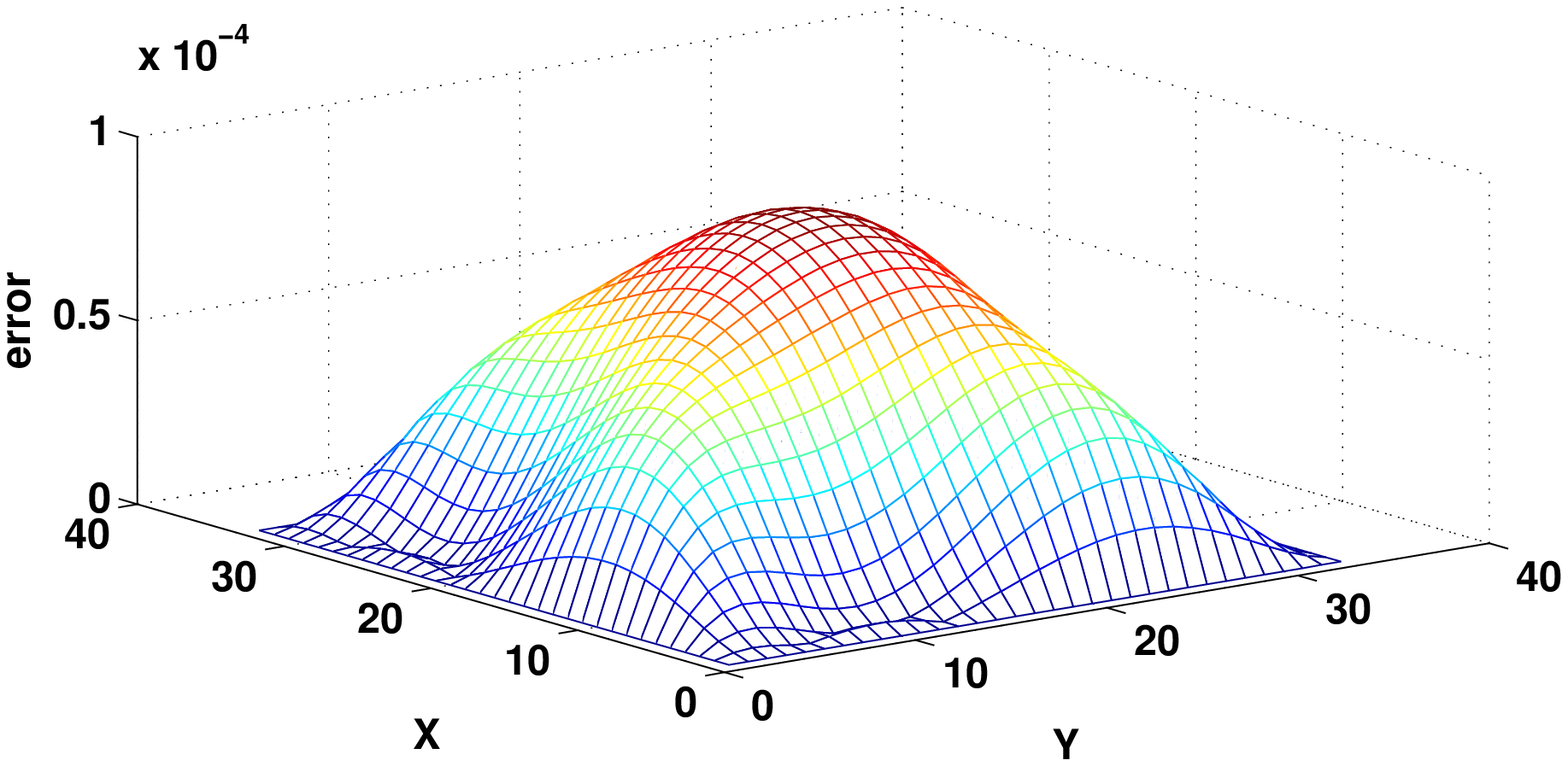,width=7.2cm}} &
\makebox{\psfig{figure=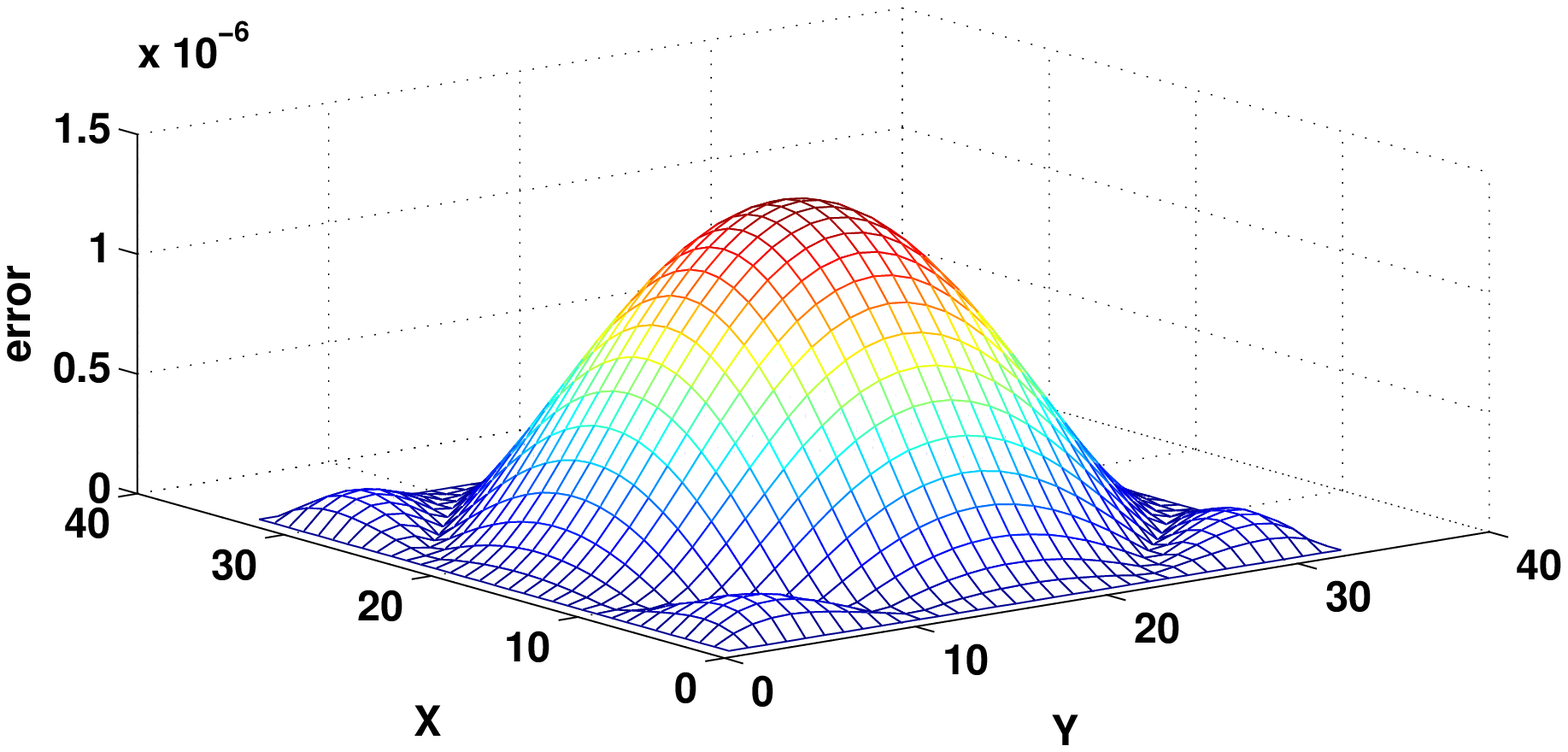,width=7.2cm}}
\\  (a) &  (b)
\end{tabular}
\caption{Error in maximum norm for the \emph{Example 1}:  (a) CDS with mesh parameters $M_x=M_y=32$,
 $N=32$; (b) CFDS  for $M_x=M_y=32$, $N=256$}\label{Fig2}
\end{figure}

In Table \ref{table2D1} the mesh refinement analysis using CDS and CFDS are presented.
The results confirm the theoretical rate of convergence, i.e. the ratio near four confirm second order for the CDS and near sixteen - fourth order for the CFDS.
Also, as the CFDS has an error $O(h^{4}+\tau^{2})$, to observe the fourth order, when doubling the number of mesh points in space  one must  take  quadruple mesh points in time.
The advantage of the CFDS is corroborated by presenting the CPU time - there needs smaller time for the CFDS to obtain results with better accuracy in despite of the using of more time layers. In Fig.~\ref{Fig1} the exact solution at final time $T$ for $u_1$ and mesh parameters $M_x=M_y=32$,
 $N=32$ is depicted. In Fig.~\ref{Fig2}  the error, obtained by a) CDS for $M_x=M_y=32$,
 $N=256$ and by b) CFDS $M_x=M_y=32$, $N=256$ are presented.
\begin{table}[tbp]
\caption{Comparison of the maximum absolute errors of the CDS and CFDS for  \textbf{Example 1}} \label{table2D1}\centering
\begin{tabular}{|c|c|c|c|c|c||c|c|c|c|c|c|}
\hline
 \multicolumn{6}{|c||}{ CDS,  $O(h^{2}+\tau^{2})$}&
 \multicolumn{6}{|c|}{CFDS,  $O(h^{4}+\tau^{2})$} \\
 \hline $M_x$ & $M_y$    & N      &   $error_M$     & $ratio$  & CPU    & $M_x$      & $M_y$    & N   & $error_M$     &  $ratio$  & CPU
\\
 \hline
4 & 4 & 4 & 5.702 e-03 &  -   & 0.58
&4 & 4 & 4 &    5.875 e-03 &  -   & 0.72
\\
8 &  8 & 8   &  1.449 e-03 &  3.94    &1.82
&8 &  8 & 16  & 3.595 e-04 & 16.34    &3.04
\\
16 & 16 & 16 &    3.637 e-04 & 3.99    &  14.42
& 16 & 16 & 64    & 2.232 e-05 & 16.11 & 29.74
\\
32 & 32  & 32  &  9.102 e-05   & 4.001   & 143.7
&  32 & 32  & 256 & 1.392 e-06   &  16.03  &1076
\\
64 & 64  & 64  & 2.276 e-05    & 4.00  & 3959
& 64 & 64  & 1024  & 8.698 e-08    & 16.003  & 60907
\\
128 & 128 & 128  & 5.691 e-06    & 4.00  & 32709
& 128 & 128  & 4096  & 5.436 e-09    & 16.0001  & 720477
\\ \hline
\end{tabular}
\end{table}

\begin{figure}[tbp]
\centering
\begin{tabular}{cc}
\makebox{\psfig{figure=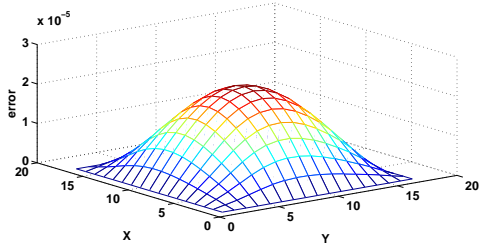,width=7.2cm}} &
\makebox{\psfig{figure=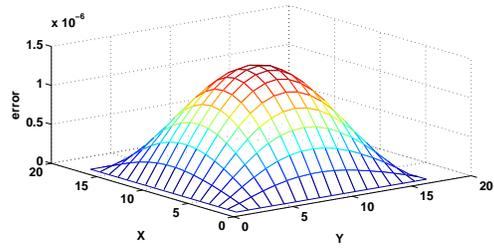,width=7.2cm}}
\\  (a) &  (b)
\end{tabular}
\caption{Error in maximum norm for the \emph{Example 1}:  (a) CDS with RE in space $M_x=M_y=16$,
 $N=64$ ; (b) CFDS with RE in space and $M_x=M_y=16$, $N=256$}\label{Fig3}
\end{figure}

In Table \ref{table2D2} the mesh refinement analysis using CDS and CFDS with Richardson extrapolation (RE) in space (using corresponding weights from (\ref{RE2}) and (\ref{RE4})) are presented. Again, to observe the fourth and sixth order of  CDSRE and CFDSRE, doubling mesh points in space one must take the number of time layers four and eight times more from the previous experiment. The results confirm the expected rates of convergence for both numerical methods. The ratio near 64 corresponds with sixth order of the CFDSRE.   Comparing of the CPU time of Table \ref{table2D1} and Table \ref{table2D2} shows a priority of using Richardson Extrapolation obtaining smaller errors for smaller computational time, nevertheless that the Richardson Extrapolation needs to compute the numerical solutions on two consecutive meshes. The advantage of CFDS  with RE is also clearly seen. In Fig.~\ref{Fig3}  the error, obtained by a) CDS with RE for $M_x=M_y=16$,
 $N=64$ and by b) CFDS  with RE in space and  $M_x=M_y=16$, $N=256$ are presented.
\begin{table}[ht]
\caption{Comparison of the errors in maximum norm for the numerical \textbf{Example 1 } for CDS and CFDS with Richardson extrapolation in space}\label{table2D2}\centering
\begin{tabular}{|c|c|c|c|c|c||c|c|c|c|c|c|}
\hline
 \multicolumn{6}{|c||}{ CDS with RE in space, $O(h^{4}+\tau^{2})$}&
 \multicolumn{6}{|c|}{CFDS with RE in space, $O(h^{6}+\tau^{2})$} \\
 \hline $M_x$ & $M_y$    & N      &   $err_N$     & ratio  & CPU    & $M_x$      & $M_y$    & N   & $err_N$     &  ratio  & CPU
\\
 \hline
4 & 4 & 4 & 5.677 e-03 &  -   & 1.34
& 4 & 4 & 4 & 5.711 e-03    &  -   & 1.38
\\
8 &  8 & 16   &   3.545 e-04 &  16.014     &16.17
& 8 &  8 & 32  &  8.912 e-05    & 64.087       & 17.45
\\
16 & 16 & 64 &    2.216 e-05 &  15.997    &  544
& 16 & 16 & 256  &  1.392 e-06     &  64.022    &  1497
\\
32 & 32  & 256  &  1.385 e-06     &  16.001     & 3055
&  32 & 32 & 2048 & 2.1757 e-08     &  63.989    & 23390
\\ \hline
\end{tabular}
\end{table}

\begin{figure}[tbp]
\centering
\begin{tabular}{cc}
\makebox{\psfig{figure=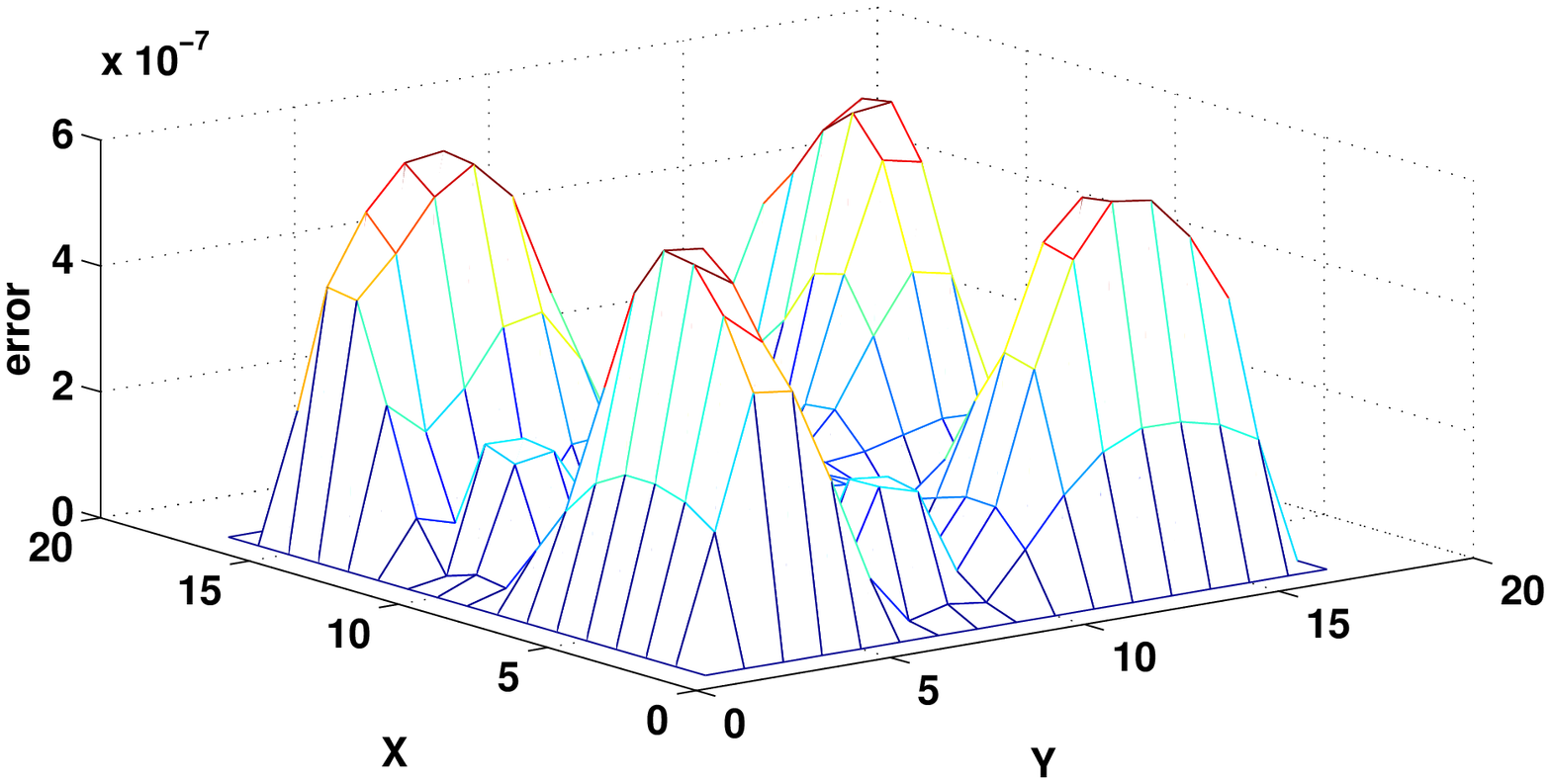,width=7.2cm}} &
\makebox{\psfig{figure=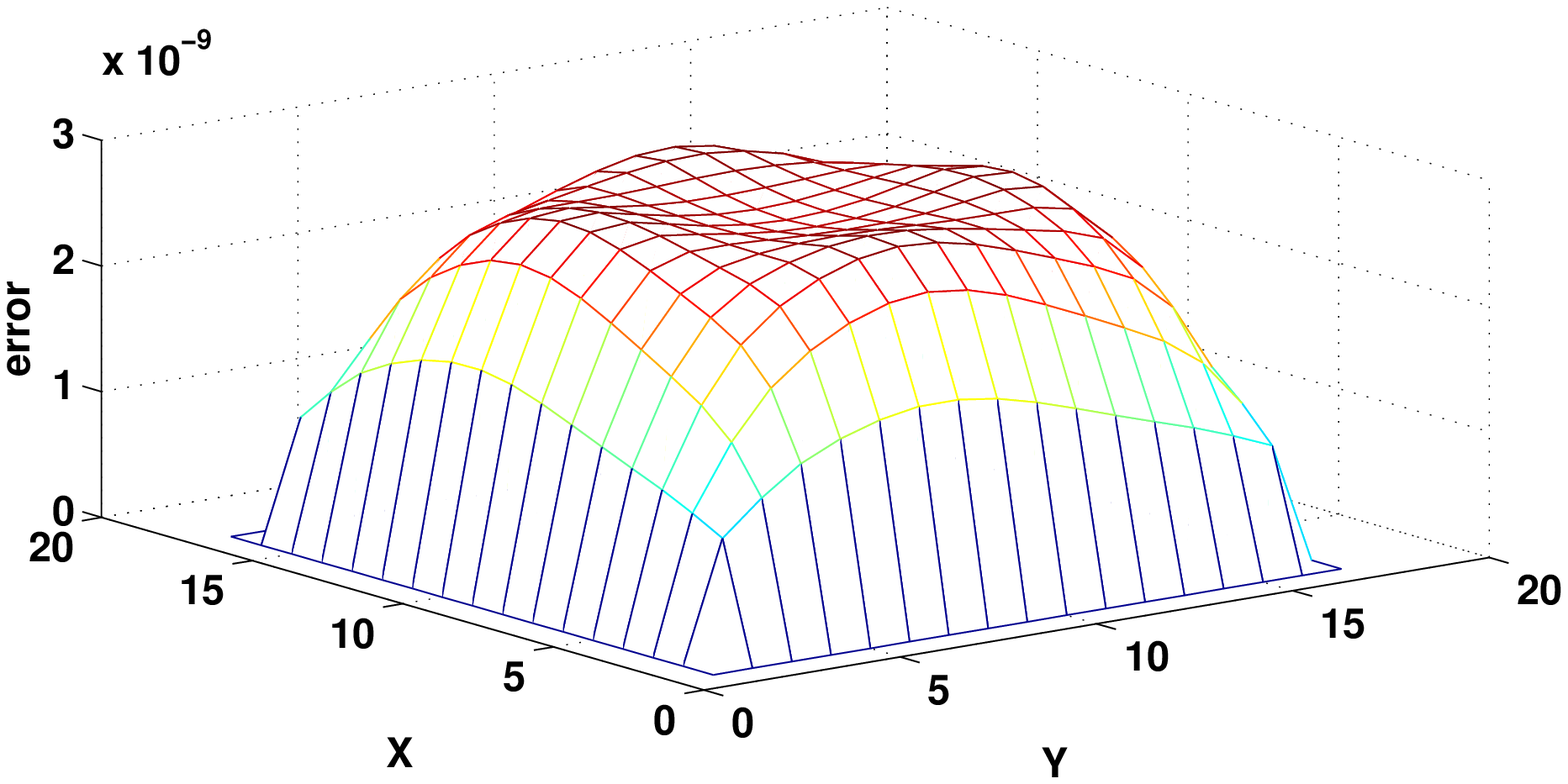,width=7.2cm}}
\\  (a) &  (b)
\end{tabular}
\caption{Error in maximum norm for \emph{Example 1}:  (a) CDS with RE in space and time $M_x=M_y=16$,
 $N=16$; (b) CFDS with RE in space and time $M_x=M_y=16$, $N=64$}\label{Fig4}
\end{figure}

In Table \ref{table2D3} the mesh refinement analyses using CDS and CFDS with Richardson extrapolation (RE) in space and time are presented. Again, to observe the fourth and sixth order of  CDSRE and CFDSRE, doubling mesh points in space one must take the number of time layers two and eight times more from the previous experiment. This would cause to extremely growth of CPU time for the case of CFDS and therefor we take here four times (instead eight times) smaller mesh intervals in time. The results confirm the expected rates of convergence for both numerical methods. Comparing of the CPU time of Tables \ref{table2D1}, \ref{table2D2} and \ref{table2D3} shows a priority of using Richardson Extrapolation both in space and time obtaining smaller errors for smaller computational time. The advantage of CFDSRE is also clearly seen. Fig.~\ref{Fig4} presents the errors in maximum norm for  \emph{Example 1}  (a) with CDS and RE in space and time $M_x=M_y=16$,  $N=16$; (b) with CFDS and RE in space and time $M_x=M_y=16$, $N=64$ and is in concurdance with the results in Table \ref{table2D3}.
\begin{table}[ht]
\caption{Comparison of the errors in maximum norm for the numerical \textbf{Example 1 } for CDS and CFDS with Richardson extrapolation in space and time }\label{table2D3}\centering
\begin{tabular}{|c|c|c|c|c|c||c|c|c|c|c|c|}
\hline
 \multicolumn{6}{|c||}{ CDS with RE in space and time, $O(h^{4}+\tau^{4})$}&
 \multicolumn{6}{|c|}{CFDS with RE in space and time $O(h^{6}+\tau^{4})$} \\
 \hline $M_x$ & $M_y$    & N      &   $err_N$     & ratio  & CPU    & $M_x$      & $M_y$    & N   & $err_N$     &  ratio  & CPU
\\
 \hline
4 & 4 & 4 & 5.649 e-05    &  -   & 6.73
& 4 & 4 & 4 &  8.476 e-06  &  -   & 3.36
\\
8 &  8 & 8  &  9.722 e-06    &  5.81     & 18.71
& 8 &  8 & 16  & 1.748 e-07& 48.49    &30.26
\\
16 & 16 & 16  &  5.989 e-07  & 16.23   & 194.81
& 16 & 16 & 64  & 2.847 e-09  & 61.39 & 1276
\\
32 & 32 & 32 &  3.715 e-08  &  16.12  & 4594
& 32  & 32  & 256  & 4.529 e-11   &  62.86   & 66991
\\
64 & 64 & 64 & 2.171 e-09   &  16.03   & 37101
&  64 & 64  & 1024  & 7.086 e-13   &  63.91   & 790800
\\ \hline
\end{tabular}
\end{table}

\subsection{Example 2 ( no exact solution)}
In this case we consider more realistic variant of problem (\ref{2Dsystem1})-(\ref{starn}) with the following parameters of the domain: the spatial domain is the square $\Omega=[0,500]^2$ with side length $500$ \emph{km}, the length of the time interval $[0,T]$ is $1440$ \emph{min} and the number of equations is $L=10$. The initial conditions on the time level $t=0$ are the constant functions
$$ \mathbf{u}_0(x,y)= (10^3,10^3,10^3, 5.10^3, 5.10^3, 10^2, 10^{-2}, 10^{-2}, 10^{-3}, 10^{-11} ),$$
measured in $mol/km^3$ and the boundary conditions are chosen to be periodic: $\gamma_i$ has the form
$$\gamma_l (t)= const_l  (sin(t/C)+2),$$
where $C=4$ is a constant and the constants $const_l$, $l=1,...,L$ are chosen in such a way that the compatibility of the boundary and initial data is ensured. The diffusion coefficient is set to be $K=1.8 km^2/min$ and the coefficient $\mu$ is $\mu=2\pi/(60*T)$.

In this example there is not analytical solution. One  way of calculating the convergence rate is the method of Runge on three nested meshes. Here we use another idea. As an "exact" solution we take the solution, obtained with a "least" mesh size in space. In the following tables we denote these solutions by \textbf{bold} font. Also in this case we present the relative error in maximum norm.
We control the rate of convergence denoted by $order$ and evaluated by
 \begin{equation*}
    order=log_2(ratio)
\end{equation*}
when doubling the number of mesh points and in other case
\begin{equation*}
    order=log(error_{M',l}/ error_{M'',l})/log(M''/M')
\end{equation*}
where $M'$ and $M''$ are two consecutive numbers of mesh points in space in the mesh refinement analysis.

In Table \ref{tableAir1} we present the results obtained by CDS with number of  time steps $N=256$ for the first and fifth substances $ u_1$ and $ u_5$ at the central node with coordinates $(x_{M/2},y_{M/2})=(X/2,Y/2)=(250,250)$.
The second order is confirmed. It is interesting to note that neverthelees $u_1$ and $u_5$ have different values, the relative errors are approximately the same for the both pollutants.
Similar results are presented in Table \ref{tableAir2}, but at the point $(x,y)=(X/6,Y/6)=(83.33,83.33)$. Again the second order of the CDS can be seen.
\begin{table}[!htb]\caption { The rate of convergence for the \textbf{Example 2} for the CDS at the central node $(x,y)=(X/2,Y/2)$  with  time steps $N=256$ for the first and fifth substances $ u_1$ and $ u_5$}\label{tableAir1}
{\tabcolsep=2pt \baselineskip=3pt 
\begin{center}
\begin{tabular}{|c|c|c|c|c||c|c|c|c|c|}
\hline
\multicolumn{5}{|c||}{ $U_1$}&
 \multicolumn{5}{|c|}{$U_5$} \\
\hline\hline $M_x$ & $M_y$ & numerical value & rel. error & order & $M_x$ & $M_y$ & numerical value & rel. error & order  \\
\hline $8$    &$8$        &1975.88248125790       &1.001 e-02         &  -     & $8$    &$8$        &4523.29297726041       &1.001 e-03         &  -     \\
\hline $16$    & $16$        &1991.14360768096       &2.366 e-03         &2.08                 & $16$    & $16$        &4558.22937850949       &2.366 e-03         &2.08                 \\
\hline $24$    &$24$        &1993.81301129742       &1.028 e-03       &2.05            & $24$    &$24$        &4564.34028098519       &1.028 e-03       &2.05            \\
\hline $32$    & $32$       &1994.73061732235       &5.685 e-04         &2.06            & $32$    & $32$       &4566.44089971161       &5.684 e-04         &2.06             \\
\hline $40$     &$40$      &1995.15232367582       &3.572 e-04         &2.08            & $40$     &$40$      &4567.40628589102       &3.572 e-04         &2.08            \\
\hline $48$    & $48$        &1995.38060726527       &2.428 e-04         &2.11           &  $48$    & $48$        &4567.92888132097       &2.428 e-04         &2.11           \\
\hline $56$    &  $56$     &1995.51798902418       &1.739 e-04         &2.16              & $56$    &  $56$     &4568.24338083461       &1.740 e-04         &2.16               \\
\hline $64$    & $64$       &1995.60704897582       & 1.293 e-04      & 2.21         &  $64$    & $64$       &4568.44726023862       & 1.293 e-04      & 2.21           \\
\hline $192$& $192$       &\textbf{1995.86518532405}           &               &                & $192$& $192$       &\textbf{4569.03819569955}           &               &                \\
\hline\hline
\end{tabular}
\end{center}
}
\end{table}

\begin{table}[!htb]\caption { The numerical values, the relative errors and the rate of convergence for \textbf{Example 2} by the CDS at the node $(x,y)=(X/6,Y/6)$  with number of  time steps $N=256$ for the first and fifth substances $ u_1$ and $ u_5$}\label{tableAir2}
{\tabcolsep=2pt \baselineskip=3pt 
\begin{center}
\begin{tabular}{|c|c|c|c|c||c|c|c|c|c|}
\hline
\multicolumn{5}{|c||}{ $U_1$}&
 \multicolumn{5}{|c|}{$U_5$} \\
\hline\hline $M_x$ & $M_y$ & numerical value & rel. error & order & $M_x$ & $M_y$ & numerical value & rel. error & order  \\
\hline $6$    &$6$        &1068.47327302014       &4.271 e-02         &  -     & $6$    &$6$        &2447.7334068223       &4.203 e-02         &  -     \\
\hline $12$    & $12$        &1110.55728440439       &5.007 e-03         &3.09                 & $12$    & $12$        &2542.36959186444       &4.998 e-03         &3.07                 \\
\hline $24$    &$24$        &1115.53721634304       &5.451 e-04       &3.19           & $24$    &$24$        &2553.74668391763       &5.450 e-04       &3.20           \\
\hline $48$    & $48$       &1116.05637283823       &7.994 e-05       &2.76             & $48$    & $48$       &2554.93507405114       &7.992 e-05       &2.77             \\
\hline $96$    & $96$       &1116.14559394767       &1.783 e-05       &2.16             & $96$    & $96$       &2555.13927953381       &1.782 e-05
&2.16             \\
\hline \hline $192$    & $192$       &\textbf{1116.16549194698}       &     -             &                & $192$    & $192$       &\textbf{2555.18481923814}       &     -             &                \\
\hline\hline
\end{tabular}
\end{center}
}
\end{table}

With the same parameters  the experiments are repeated using CFDS. The results are presented in Table \ref{tableAir3} and Table \ref{tableAir4}. The fourth order in both cases (central node (x,y)=(X/2,Y/2) and node (x,y)=(X/6,Y/6) )  for the both substances  $u_1$ and $u_5$ is confirmed. Again at the central node the relative errors are likely the same.

\begin{table}[tbp]\caption {The rate of convergence for \textbf{Example 2} for the CFDS at the central node $(x,y)=(X/2,Y/2)$  with number of time steps $N=256$ for the first and fifth substances $ u_1$ and $ u_5$}\label{tableAir3}
{\tabcolsep=2pt \baselineskip=3pt 
\begin{center}
\begin{tabular}{|c|c|c|c|c||c|c|c|c|c|}
\hline
\multicolumn{5}{|c||}{ $U_1$}&
 \multicolumn{5}{|c|}{$U_5$} \\
\hline\hline $M_x$ & $M_y$ & numerical value & rel. error & order & $M_x$ & $M_y$ & numerical value & rel. error & order  \\
\hline $8$    &$8$        &2000.63329684645       &2.273 e-03     &  -    & $8$    &$8$        &4580.15403342582       &2.417 e-03         &  -     \\
\hline $16$    & $16$        &1996.19582729555       &1.495 e-04     &3.988       & $16$& $16$     &4569.79512223399       &1.495 e-04         &4.014                 \\
\hline $24$    &$24$        &1995.95673693047       &2.972 e-05       &3.984           & $24$    &$24$        &4569.24777942408       &2.972 e-05       &3.984            \\
\hline $32$    & $32$       &1995.91621977218       &9.419 e-06         &3.994             & $32$    & $32$       &4569.15502567325       &9.420 e-06         &3.994              \\
\hline $40$     &$40$      &1995.90511881075       &3.858 e-06         &4.000            &  $40$     &$40$     &4569.12961288337       &3.858 e-06         &4.000             \\
\hline $48$    & $48$        &1995.90112653914       &1.858 e-06         &4.008               & $48$    & $48$       &4569.12047362402       &1.858 e-06         &4.008               \\
\hline $56$    &  $56$     &1995.89941400964       &9.997 e-07         &4.019               & $56$    &  $56$    &4569.11655324515       &9.997 e-07         &4.019                   \\
\hline $64$    & $64$       &1995.89858249863       &5.831 e-07       &4.037              & $64$    & $64$       &4569.11464972741       &5.831 e-07       &4.037              \\
\hline $192$    & $192$       &\textbf{1995.89741860066}       &               &                & $192$    & $192$       &\textbf{4569.11212069383}       &               &                \\
\hline\hline
\end{tabular}
\end{center}
}
\end{table}
\begin{table}[!htb]\caption { The rate of convergence of the CFDS for \textbf{Example 2} at the  node $(x,y)=(X/6,Y/6)$  with number of time steps $N=256$ for the first and fifth substances $u_1$ and $ u_5$}\label{tableAir4}
{\tabcolsep=2pt \baselineskip=3pt 
\begin{center}
\begin{tabular}{|c|c|c|c|c||c|c|c|c|c|}
\hline
\multicolumn{5}{|c||}{ $U_1$}&
 \multicolumn{5}{|c|}{$U_5$} \\
\hline\hline $M_x$ & $M_y$ & numerical value & rel. error & order & $M_1$ & $M_y$ & numerical value & rel. error & order  \\
\hline $6$    &$6$        &1043.29329103805       &6.529 e-02         &  -     & $6$    &$6$        &2257.94831662249       &1.163 e-01         &  -     \\
\hline $12$    & $12$        &1118.08045908966       &1.710 e-03         &5.25                 & $12$    & $12$        &2550.11225999479       &1.991 e-03         &5.86                 \\
\hline $24$    &$24$        &1116.07801239889       &8.411 e-05       &4.34          & $24$    &$24$        &2554.99929160301           &7.834 e-05       &4.66           \\
\hline $48$    & $48$       &1116.16605487952       &5.229 e-06       &4.00             & $48$    & $48$       &2555.18608471763       &5.236 e-06       &3.91            \\
\hline $96$    & $96$       &1116.17155052937       &3.054 e-07         &4.09                & $96$    & $96$       &2555.19868000388       &     3.068 e-07        &    4.09           \\
\hline $192$& $192$       &\textbf{1116.17189141636}       &     -             &                & $192$    & $192$       &\textbf{2555.19946394382 }      &     -             &                \\
\hline\hline
\end{tabular}
\end{center}
}
\end{table}

In Table \ref{tableAir5}  and Table \ref{tableAir6}  the results obtained by the CDSRE and CFDSRE in space are shown. The number of time layers are $N=256$ and the presented values are the numerical values at the last time layer $t_N=T$ at the central node $(x,y)=(X/2,Y/2)$. The results confirm the fourth order for the CDSRE and sixth order for CFDSRE.

\begin{table}[!htb]\caption { The rate of convergence for \textbf{Example 2} for the CDS with RE in space at the central node $(x,y)=(X/2,Y/2)$ with time steps $N=256$}\label{tableAir5}
{\tabcolsep=2pt \baselineskip=3pt 
\begin{center}
\begin{tabular}{|c|c|c|c|c||c|c|c|c|c|}
\hline
 \multicolumn{5}{|c||}{ $U_1$}&
 \multicolumn{5}{|c|}{$U_5$} \\
\hline\hline $M_x$ & $M_Y$ & numerical value & rel. error & order & $M_x$ & $M_Y$ & numerical value & rel. error & order  \\
\hline $8$    &$8$        &1996.23064982198   &1.669 e-04         &  -    & $8$    &$8$  &4569.87484559225 &1.669 e-04      &  -     \\
\hline $16$    & $16$        &1995.92628720281   &1.446 e-05         &3.529  & $16$    &$16$ &4569.17807344565 &1.446 e-05  &3.529    \\
\hline $24$    &$24$        &1995.90313925455   &2.862 e-06           &3.995  & $24$    &$24$ &4569.12508143290 &2.863 e-06  &3.994  \\
\hline $32$    & $32$       &1995.89919286031   &8.856 e-07         &4.079  & $32$    &$32$ &4569.11604708096 &8.855 e-07    &4.078 \\
\hline $40$     &$40$     &1995.89812877923   &3.524 e-07         &4.129  & $40$    &$40$ &4569.11361111869 &3.524 e-07    &4.129  \\
\hline $48$    & $48$       &1995.89775101360   &1.632 e-07          &4.225  & $48$    &$48$ &4569.11274631368 &1.631 e-07 &4.224   \\
\hline $56$    & $56$       &1995.89759041518   &8.268 e-08          &4.409  & $56$    &$56$ &4569.11237866122 &8.267 e-08 &4.409   \\
\hline $64$    & $64$       &1995.89751293823   &4.386 e-08       &4.748  & $64$& $64$ &4569.11220125629    &4.385 e-08    &4.749     \\
\hline $96$    & $96$       &\textbf{1995.89742540660}      &          &      & $96$    & $96$      &\textbf{4569.11200091091 }      &         &     \\
\hline\hline
\end{tabular}
\end{center}
}
\end{table}

\begin{table}[!htb]\caption {The rate of convergence for \textbf{Example 2} for CFDS with RE in space for the central node $(x,y)=(X/2,Y/2)$ with time steps $N=256$ }\label{tableAir6}
{\tabcolsep=2pt \baselineskip=3pt 
\begin{center}
\begin{tabular}{|c|c|c|c|c||c|c|c|c|c|}
\hline
 \multicolumn{5}{|c||}{ $U_1$}&
 \multicolumn{5}{|c|}{$U_5$} \\
\hline\hline $M_x$ & $M_y$ & numerical value & rel. error & rate & $M_x$ & $M_Y$ & numerical value & rel. error & rate  \\
\hline $8$    &$8$        &1995.89999599216     &1.299 e-06   &  -   & $8$    &$8$        &4569.10452815453       &1.624 e-06         &  -   \\
\hline $16$    & $16$        &1995.89757927062     &8.779 e-08   &3.887 & $16$& $16$        &4569.11235256920       &8.767 e-08         &4.212     \\
\hline $24$    &$24$        &1995.89741917972     &7.582 e-09   &6.040 & $24$& $24$     &4569.11198657068       &7.565 e-09       &6.042  \\
\hline $32$    & $32$       &1995.89740668039     &1.300 e-09   &6.077 & $32$& $32$       &4569.11195799769       &1.311 e-09         &6.092  \\
\hline $40$     &$40$     &1995.89740473047     &3.428 e-10   &6.041 & $40$&  $40$    &4569.11195355115      &3.380 e-10         &6.075 \\
\hline $48$    & $48$       &1995.89740427467     &1.144 e-10   &6.018 & $48$& $48$       &4569.11195251649       &1.115 e-10         &6.080     \\
\hline $56$    & $56$       &1995.89740413622     &4.505 e-11   &6.045 & $56$& $56$       &4569.11195220449       &4.326 e-11         &6.144     \\
\hline $64$    & $64$       &1995.89740408556     &1.968 e-11   &6.204 & $64$& $64$       &4569.11195209159       &1.855 e-11         &6.339     \\
\hline \hline$96$    & $96$       &\textbf{1995.89740404629}      &          &      & $96$    & $96$       &\textbf{4569.11195200681}       &               &     \\
\hline\hline
\end{tabular}
\end{center}
}
\end{table}

\begin{figure}[htbp]
\centering
\psfig{figure=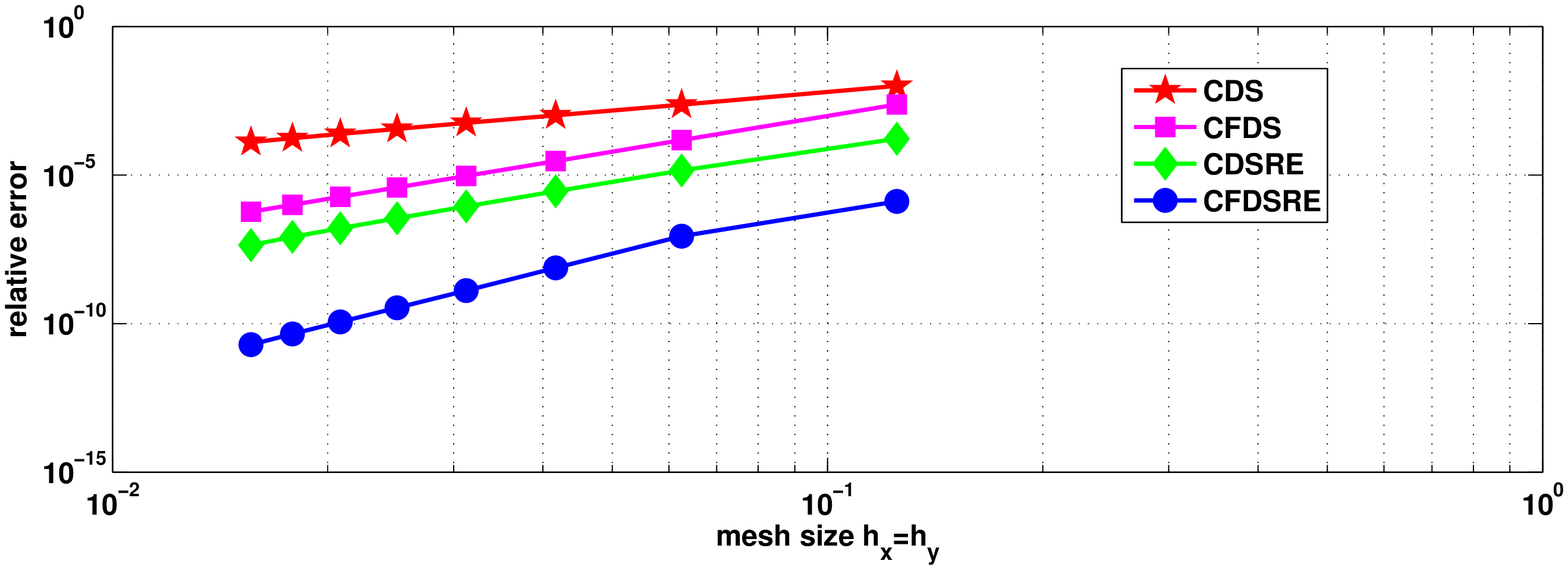,width=10.0cm}
\caption{The log-log plot of the errors versus space mesh size for the \emph{Example 2}, obtained by:  CDS - {\color{red} red line}, ${\color{red}-\star-}$;  CFDS - {\color{magenta}magenta line}, ${\color{magenta}-\blacksquare-}$; CDSRE - {\color{green} green line}, ${\color{green}-\blacklozenge-}$; CFDSRE in space -  {\color{blue} blue line}, ${\color{blue} -\bullet-}$.}\label{loglog}
\end{figure}

In Fig.~\ref{loglog} the log-log plot of the errors versus space mesh size for the \emph{Example 2} is presented, obtained by:  CDS - red line, ${\color{red}-\star-}$;  CFDS - magenta line, ${\color{magenta}-\blacksquare-}$; CDSRE - green line, ${\color{green}-\blacklozenge-}$; CFDSRE in space - blue line, ${\color{blue} -\bullet-}$. The increasing of the slope of the lines corresponds with the increasing of the rate of convergence. The lowest  line confirms the advantage of the CFDS in combination with Richardson extrapolation.

\begin{figure}[tbp]
\centering
\begin{tabular}{cc}
\makebox{\psfig{figure=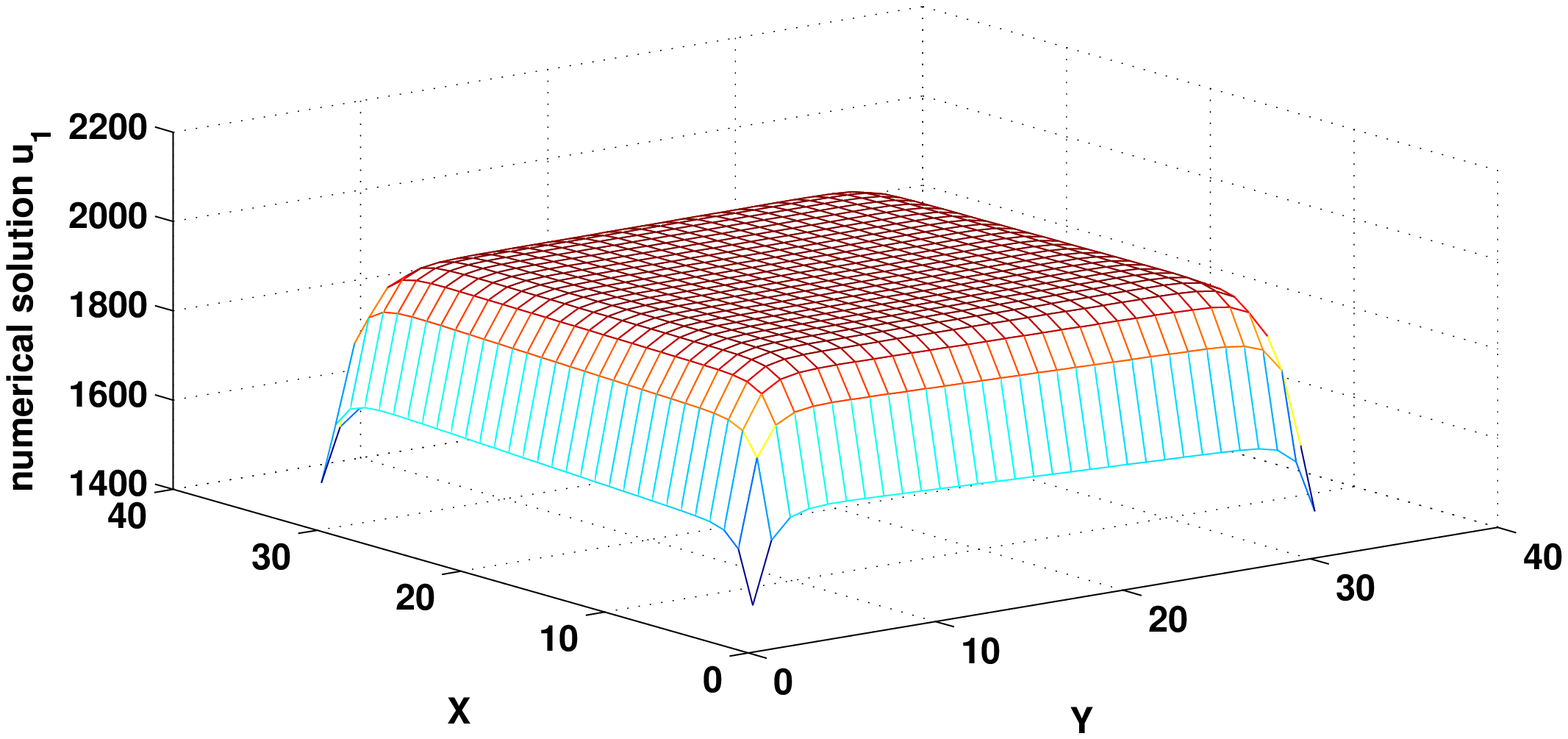,width=7.2cm}} &
\makebox{\psfig{figure=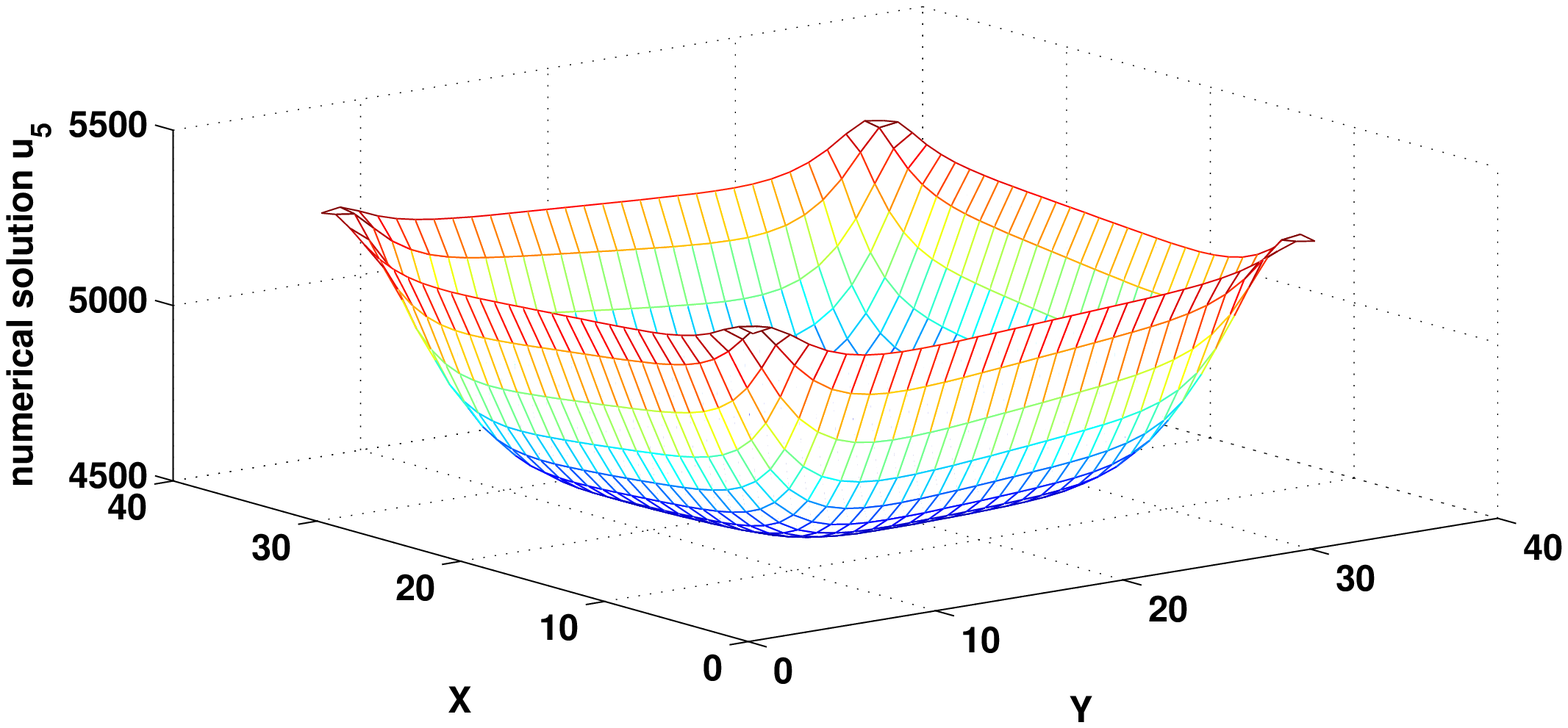,width=7.2cm}}
\\  (a) &  (b)
\end{tabular}
\caption{Numerical solution for \emph{Example 2}, obtained by CDS for $\mu=2\pi/(60T)$ with mesh parameters $M_y=M_y=32$,
 $N=256$:  (a) for $u_1$; (b) for $u_5$}\label{Fig1n}
\end{figure}

\begin{figure}[tbp]
\centering
\begin{tabular}{cc}
\makebox{\psfig{figure=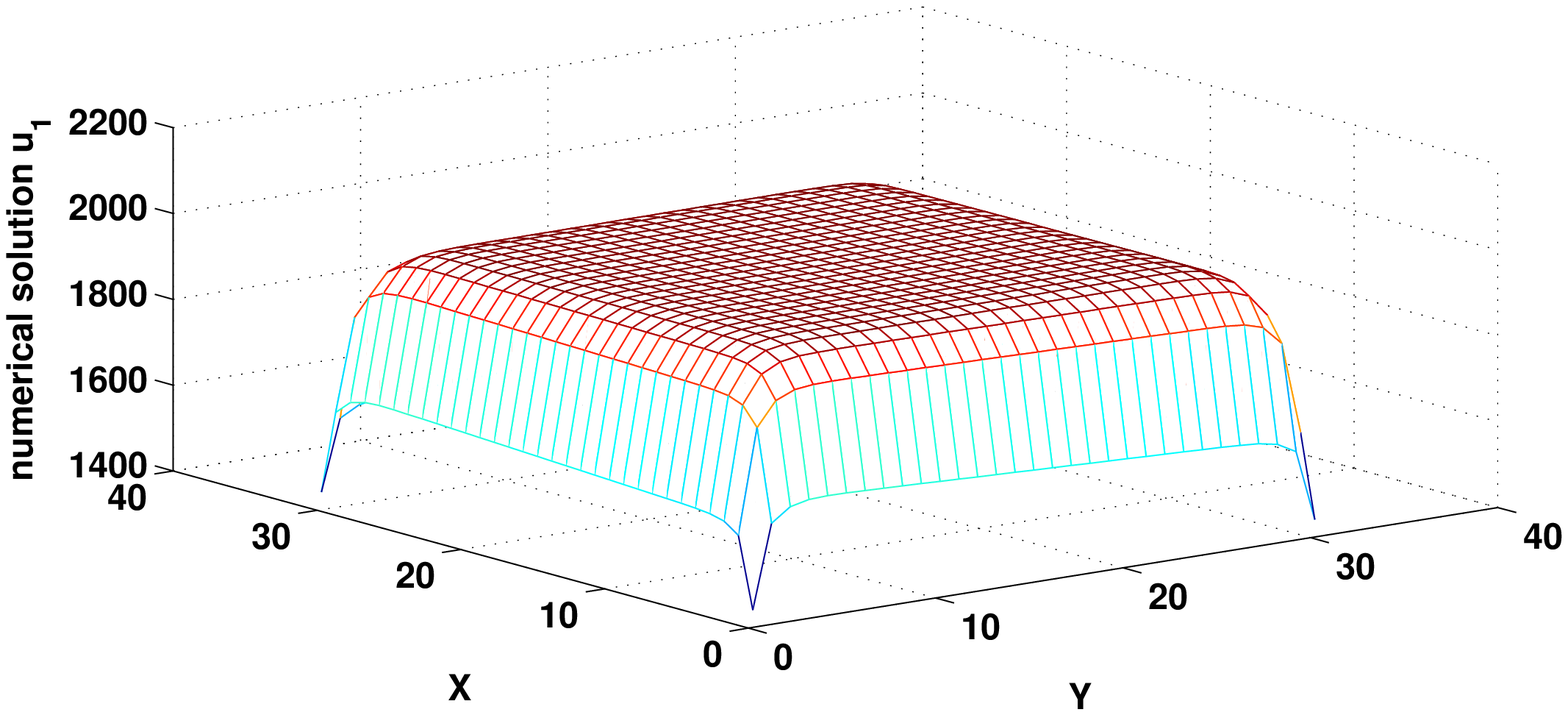,width=7.2cm}} &
\makebox{\psfig{figure=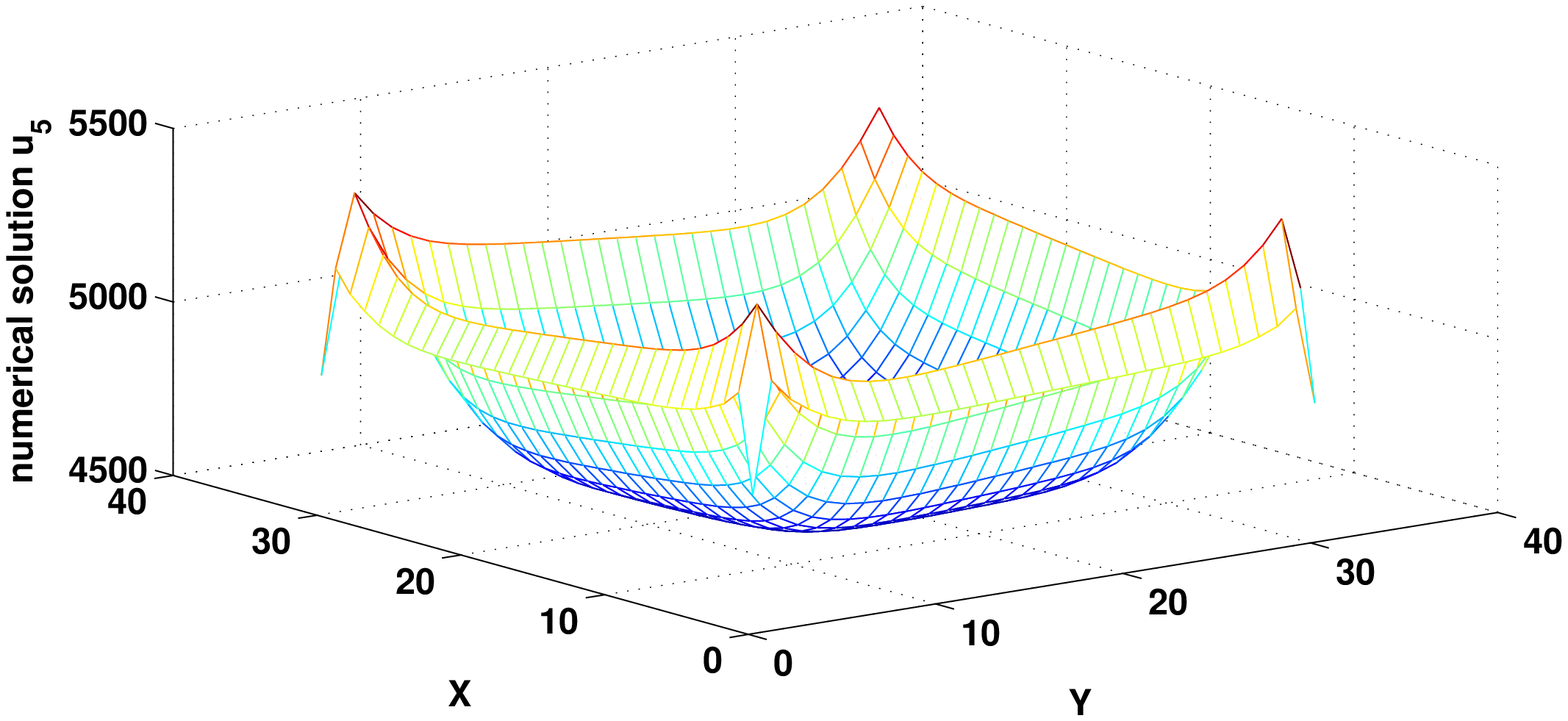,width=7.2cm}}
\\  (a) &  (b)
\end{tabular}
\caption{Numerical solution obtained by CFDS for $\mu=2\pi/(60T)$  with mesh parameters $M_x=M_y=32$,
 $N=256$ \emph{Example 2}:  (a) for $u_1$; (b) for $u_5$}\label{Fig2n}
\end{figure}

In Fig.~\ref{Fig1n} the numerical solutions obtained by CDS for $\mu=2\pi/(60T)$ with mesh parameters $M_x=M_y=32$,
at final time layer $N=256$  (a) for $u_1$; (b) for $u_5$ are shown. Similarly, in Fig.~\ref{Fig2n} the numerical solutions obtained by CFDS are shown.

\begin{figure}[htbp]
\centering
\begin{tabular}{cc}
\makebox{\psfig{figure=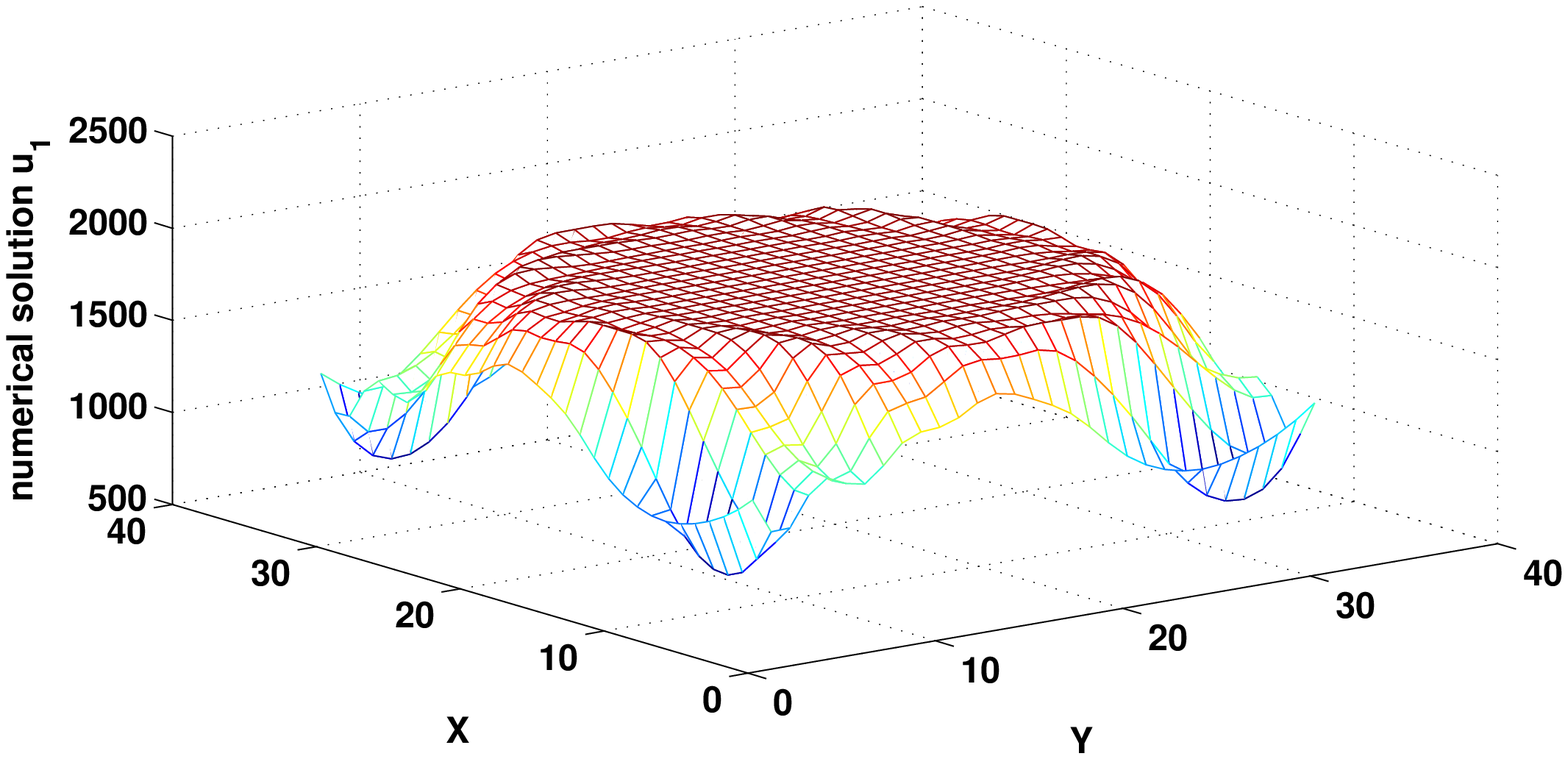,width=6.5cm}} &
\makebox{\psfig{figure=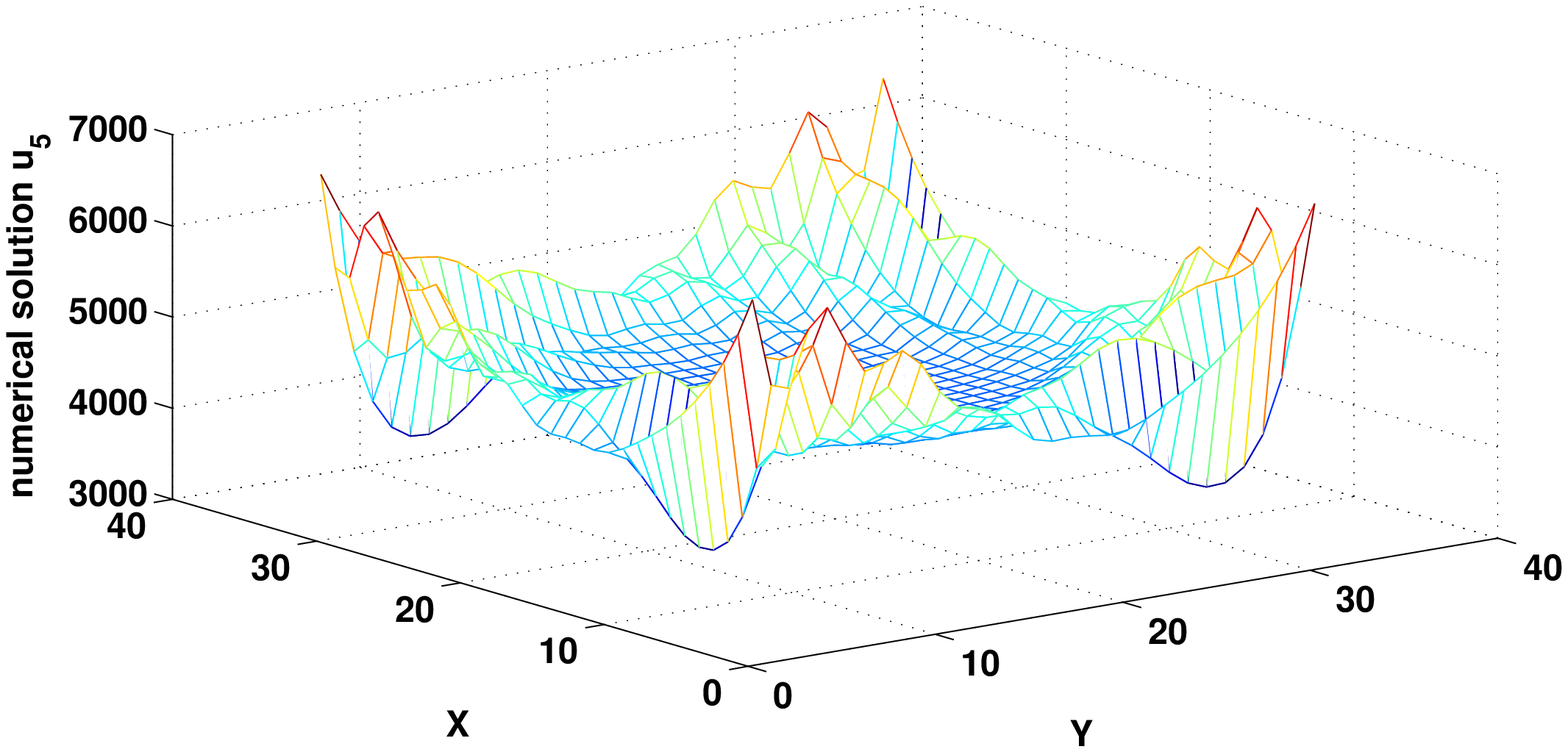,width=6.5cm}}
\\  (a) &  (b)
\end{tabular}
\caption{Numerical solution for \emph{Example 2}, obtained by CFDS with $\mu=2\pi/(500)$ and $M_x=M_y=32$,
 $N=256$:  (a) for $u_1$; (b) for $u_5$}\label{Fig2pi500}
\end{figure}

Many others experiments have been done. It is interesting to see the behaviour  of the solutions if the coefficient $\mu$ in the convection term is taken to be $\mu=2\pi/(X)$ as it is in \cite{KarKur} instead $\mu=2\pi/(60*T)$ as it is in \cite{GeoZla}.  The increasing of the convective coefficients leads to significant change of the numerical solution near the corners, see Fig.~\ref{Fig2pi500} where $\mu=2\pi/500$. It can be seen that the constant initial values have been left relatively intact in the middle of the domain, but they have been stretched near the boundary by the sinusoidal boundary conditions.

In Table~\ref{bicgnewton1} the average number of iterations for  \textit{Example 1} at the outer (Newton) and at the inner (bicgstabl) part of the inexact Newton method for CDS and CFDS are presented. To go from the $n$-th time layer to the next $n+1$-th time layer we need of approximately three iterations at the outer (Newton) part for the both difference schemes. At the inner (bicgstabl) part for the case of  CDS  we need of  three iterations and for the case of  CFDS we observe the decreasing of the number of iterations from $3.40$ to $2.05$ when the numbers of the mesh points in space and time are increasing. Similar results are presented in Table~\ref{bcignewton2} for \textit{Example 2} obtained  with the number of time steps $N=256$. The number of the outer iterations is  three for CDS and decreases from $3.80$ to $3.17$ for CFDS.
In the opposite the number of the inner (bicgstabl) iterations increases for CDS from $1.75$ to $6.54$ and decreases from $4.70$ to $2.50$ for CFDS as a result of better local approximation.
\begin{table}[tbp]
\caption{The average number of iterations for  \textit{Example 1} at the outer (Newton) and inner (bicgstabl) parts of the inexact Newton method for CDS and CFDS} \label{bicgnewton1}\centering
\begin{tabular}{|c|c|c|c|c||c|c|c|c|c|}
\hline
 \multicolumn{5}{|c||}{ CDS }&
 \multicolumn{5}{|c|}{CFDS   } \\
 \hline $M_x$ & $M_y$    & N      &   Newton     & bicgstabl      & $M_x$      & $M_y$    & N   & Newton     &  bicgstabl \\
 \hline
8 &  8 & 8   &  3 &  2.67
&8 &  8 & 16  & 3 & 3.40
\\
16 & 16 & 16 &    3  & 2.67
& 16 & 16 & 64    & 2.98 & 2.57
\\
32 & 32  & 32  &  3   & 2.67
&  32 & 32  & 256 & 2.96   & 2.15
\\
64 & 64  & 64  &  2.95   & 3.31
& 64 & 64  & 1024  &  2.65   & 2.05
\\ \hline
\end{tabular}
\end{table}
\begin{table}[!htb]\caption { The number of average iterations for  \textit{Example 2} on the outer (Newton) and inner (bicgstabl) part of the inexact Newton method for CDS and CFDS with the number of time steps $N=256$}\label{bcignewton2}
{\tabcolsep=2pt \baselineskip=3pt 
\begin{center}
\begin{tabular}{|c|c|c|c||c|c|c|c|}
\hline
\multicolumn{4}{|c||}{ $CDS$}&
 \multicolumn{4}{|c|}{$CFDS$} \\
\hline\hline $M_x$ & $M_y$ & Newton & bicgstabl  & $M_x$ & $M_y$ & Newton & bicgstabl  \\
\hline $8$    &$8$        &3       &1.75          & $8$    &$8$        &3.80   &4.70 \\
\hline $16$    & $16$        &3       &2.48          & $16$    & $16$    &3.96 &4.36  \\
\hline $32$    &$32$        &3      &3.86           & $32$    &$32$    &3.32  &3.67  \\
\hline $64$    & $64$       & 3      &  6.54        & $64$  & $64$   & 3.17 & 2.50 \\
\hline\hline
\end{tabular}
\end{center}
}
\end{table}

\begin{figure}[tbp]\label{FigNonPos}
\centering
\begin{tabular}{cc}
\makebox{\psfig{figure=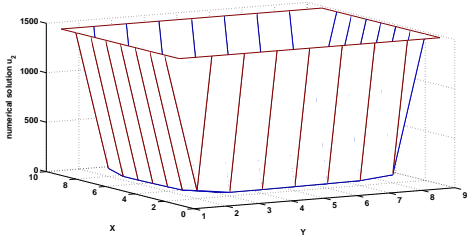,width=7.0cm}} &
\makebox{\psfig{figure=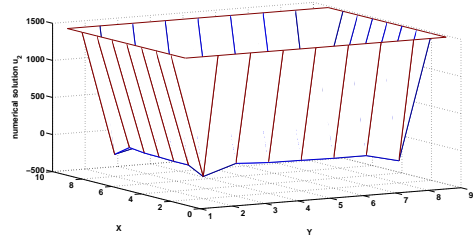,width=7.0cm}}
\\  (a) &  (b)
\end{tabular}
\caption{Numerical solution for the pollutant $NO_2$ - $u_2$  for $\mu=2\pi/(T)$  with mesh parameters $M_x=M_y=8$,
 $N=256$ \emph{Example 2}, obtained by:  (a) CDS; (b) CFDS }
\end{figure}

In spite of all advantages of CFDS in sense of accuracy and CPU time, there is also some disadvantages. The stencil of the CFDS is nine-point and the sign condition of the discrete maximum principle is not fulfill. As a result the positivity of the numerical solution is break for some values of the mesh parameters in space ant time. In Fig. \ref{FigNonPos} the numerical solution for the pollutant $NO_2$  ($u_2$)  for \emph{Example 2} when $\mu=2\pi/(T)$  and $M_x=M_y=8$, $N=256$, obtained by (a) CDS and by (b) CFDS is presented.
The CDS preserves the positivity of the numerical solution, while the CFDS does not - near the corners the numerical solution is negative and has no chemical sense. This fact confirm, that the proposed methods needs of more careful analysis.

\section{CONCLUDING REMARKS}

In this article, two different ways for derivation high-order difference schemes for semilinear parabolic systems of equations are analyzed. First, using central difference  approximation with Richardson extrapolation a fourth-order method is derived. Second, a recent proposed fourth-order in space compact difference scheme \cite{DKTVHis} is extended by Richardson extrapolation to sixth-order approximation. The time-stepping is realized using $\theta$-scheme, but in the numerical computations - by the Crank-Nicolson/Newton algorithm. The reported computational results demonstrate that the convergence rate of the CDS is  $O(h^{2}+\tau^{2})$ and of the CFDS  it is  $O(h^{4}+\tau^{2})$, but in combination with Richardson extrapolation they are respectively $O(h^{4}+\tau^{2})$ and $O(h^{6}+\tau^{2})$.
Numerically it is confirm the advantages of the CFDS over the CDS both in the accuracy and CPU time. The skilfully application of Richardson extrapolation also plays important role in obtaining good results in real time with a small number of grid nodes despite the large intervals of the domain both in space and time in air pollution problems.

In the next study we will present a theoretical analysis of the present approximation. Also, we will develop two-grid algorithms for solution of the corresponding nonlinear systems of algebraic equation. In our future work we will exploit this strategy of combining fourth-order compact difference scheme with Richardson extrapolation for solving steady-state nonlinear problems.

{\bf{Acknowledgement}}.
This work was partially supported by the Bulgarian National Fund of Science under the grant
DFNI I02-20/2014, as well as by the Program for career development of the Young scientists, BAS,
Grant No. DFNP-91/04.05.2016 and by the Project 2016-FNSE-03 of the University of Ruse.


\begin{thebibliography}{99}

\bibitem{CheLi} Chen, W., Li, C., Wright E.: On a nonlinear parabolic system-modeling
chemical reactions in rivers, Comm. on Pure and Appl.  Anal., \textbf{4}(4) 889–-899 (2005)


\bibitem{CheKin}Cheney W. , D. Kincard, \emph{Numerical Mathematics and Computing}, 4th Ed., Brooks/Cole Publishing,
Pacific Grove, CA, 1999

\bibitem{DemSta} Dembo R. S. , S. C. Eisenstat, T. Steihaug, Inexact Newton methods.  SIAM Journal on Numerical Analysis,
\textbf{19}(2), 400-408 (1982)

\bibitem{DKTVHis} Dimov, I. , Kandilarov, J., Todorov, V., Vulkov, L.:
Analysis and realization of compact difference schemes for semilinear parabolic systems, in
\emph{Numerical Methods for Scientific Computations and Advanced Applications} (NMSCAA’16), edited by K. Georgiev (Fastumprint, Sofia) 17--20 (2016)

\bibitem{DimZla} Dimov, I. , Zlatev, Z.:  \emph{Computational and Numerical Challenges
 in Air Polution Modelling.} Elsevier Science, Amsterdam-Boston-...-Tokyo, (2006).

\bibitem{GeoZla} Georgiev, K., Zlatev, Z.: Implementation of sparse matrix algorithms
in an advection-diffusion-chemistry model. J. of Comp. Appl. Math.,
\textbf{236} (3), 342-353 (2011)

\bibitem{GupManSte}
Gupta, M. M., Manohar, R. P., Stephenson, J. W.: A single cell high order scheme for the convection-diffusion equation with variable coefficients,
Int. J. for Num. Methods in Fluids,  \textbf{4}, 641-651 (1984)

\bibitem{GusKreOli} Gustafsson, B., Kreiss, H., Oliger, J.: Time Dependent Problems and Difference Methods,
Wiley, New York (1995)

\bibitem{KarKur} Karatson, J., Kurics. T.:  A preconditioned iterative solution scheme
for nonlinear parabolic systems arizing in air pollution modeling.
Math. Modell. Anal. \textbf{18} (5),  641-653 (2013)

\bibitem{KeyRoo} Kyei, Y., Roop, J. P. ,  Tang, G.:
A family of sixth-order compact finite-
difference schemes for the three-dimensional Poisson equation.
Advances in Numerical Analysis, \textbf{2010}, 1-17 (2010).

\bibitem{MarSha}Marchuk, G. I., V. V. Shaidurov,  \emph{Difference Methods and Their Extrapolations} (Springer-Verlag, New York Inc. 1983) .

\bibitem{Pao}Pao,  C. V.:    Nonlinear Parabolic and Elliptic Equations, Springer, US (1992).

\bibitem{Rich} Richards, S.: Completed Richardson extrapolation in space and time. Commun. Numer. Meth. Engineering \textbf{13}, 573--582 (1997).

\bibitem{RichtMor} Richtmyer, R.D., Morton, K.W. :  \emph{Difference Methods for
Intial-Value Problems} (Krieger, Malabar, FL., 1994).

\bibitem{Sam} A. A. Samarskii, \emph{The Theory of Difference Schemes} (Marcel Dekker, Inc.
New York, NY 2001).



\bibitem{SpoGartime} Spotz, W.,  Carey, G. F.: Extension of high-order compact schemes to time-dependent problems.
 Numer. Meth. PDE   \textbf{17}(6), 657–672 (2001)

\bibitem{SunZha}H. Sun and J. Zhang, A high-order finite difference discretization strategy based on
extrapolation for convection diffusion equations, Numer. Meth. PDEs \textbf{20}, 18–32 (2004).

\bibitem{TurGor} E. Turkel, D. Gordon, R. Gordon, and S. Tsynkov,
Compact 2D and 3D sixth
order schemes for the Helmholtz equation with variable wave number.
J. Comput. Phys. \textbf{232}, 272-287 (2013).

\bibitem{Wang} Y. Wang, High accuracy multiscale multigrid computation for partial differential equations,
Ph.D. thesis, University of Kentucky, Lexington, KY, 2010.


\bibitem{WangGuoWu} Wang, Y.-M., B.-Y. Guo, Wu, W.-J.: Fourth-order
compact finite difference methods and monotone iterative algorithms
for semilinear elliptic boundary value problems, Computers and Math.
with Appl.,\textbf{ 68}, 1671-1688  (2014)


\bibitem{ZlaDim} Zlatev, Z., Dimov, I., Farago, I., Georgiev, K. and Havasi. A.:
Application of Richardson extrapolation for multi-dimensional
advection equations, Comp. Math. Appl., \textbf{67}, 2279-2293 (2014)

\end{thebibliography}
\end{document}